%% file: integratordocu.tex
\documentclass[a4paper,12pt,twoside]{gjcommon}

\usepackage[left=2cm,right=2cm,top=2cm,bottom=3cm,includeheadfoot,bindingoffset=0.5cm]{geometry}

\renewcommand{\thesection}{\arabic{section}}

\usepackage{matlabcodehelper}
\usepackage{matlabcode}
\usepackage{booktabs}

\usepackage{hyperref}
\usepackage{setspace}
\usepackage[T1]{fontenc}
\usepackage[utf8]{inputenc}

\input{commonDefs}

\Autor{Georg Jansing}{jansing@am.uni-duesseldorf.de}
\Typ{Documentation}
\Titel{\EXPODE{} - Advanced Exponential Time Integration Toolbox for \MATLAB}
\title{\EXPODE{} - Advanced Exponential Time Integration Toolbox for \MATLAB}
\pdfsetup

\begin{document}
	\maketitle
	\renewcommand{\contentsname}{Table of Contents}
	\tableofcontents
% 	\addcontentsline{toc}{section}{Table of Contents}
	\pdfbookmark[1]{Table of Contents}{toc}
	
	\Abschnitt{Intoduction}{introduction}
	
	This document contains the documentation for the \MATLAB toolbox \EXPODE. The toolbox provides advanced exponential integration methods featuring five different integrator classes, evaluation of the Matrix functions directly and by a Krylov subspace method and an adaptive step size implementation for some of the integrators. \EXPODE is as compatible as possible to \MATLAB's internal ODE toolbox syntactically and semantically, such that a switch to this package is made easy.
	
	For the mathematical details on the implementation we refer to \cite{HocJansing10TOMS}. We first give a short quick start guide in the remainder of this section and continue with the full documentation in the next section. Note that some information might be duplicated.
	%We first give a short quick start guide and then continue with a full documentation.
	
	\input{install-quickstart}

	\input{docutext}

	\bibliographystyle{plain}
	\bibliography{integratordocu}
	\addcontentsline{toc}{section}{Literature}
\end{document}

%% file: commonDefs.tex
\newcommand{\code}[1]{\texttt{#1}}
\newcommand{\flag}[1]{\str{#1}\xspace}

\newcommand{\file}[1]{\code{#1}}
\newcommand{\newpackagename}[1]{%
	\expandafter\def\csname #1\endcsname{\code{#1}\xspace}%
}
\newpackagename{MATLAB}
\newpackagename{OCTAVE}
\newpackagename{BASH}
\newpackagename{EXPODE}
\newpackagename{EXPODEDEV}

\newcommand{\EXPODECMD}{\texttt{expode}\xspace}

\newcommand{\newmethod}[1]{%
	\expandafter\def\csname #1\endcsname{\code{\lowercase{#1}}\xspace}%
}
\newcommand{\EXP}[1]{\texttt{exp#1}\xspace}
\newmethod{EXPRK}%
\newmethod{EXPRB}%
\newmethod{EXPMS}%
\newmethod{EXPMSSEMI}%
\newmethod{EXPNEW}%
\newmethod{ODE}%
\newmethod{ODESET}%
\newmethod{EXPSET}%
\newmethod{EXPRBSET}%
\newmethod{EXPRKSET}%
\newmethod{OPTINFO}%
\newmethod{EXPRBINFO}%
\newmethod{DEVALEXP}%
\newmethod{DEVAL}%

%% file: install-quickstart.tex
\subsection{Installation and Requirements}

We will now describe the minimal requirements and the installation of the \EXPODE toolbox. \EXPODE runs on all recent and middle-aged computers. The performance highly depends on the problem and available hardware. The toolbox was tested on \MATLAB versions down to \MATLAB 7.2 (R2006a), released in 2006. Older versions might be albe to run it as well and we would appreciate feedback on success or failure. Versions prior to 7.0 cannot be compatible due to the lack of proper function handles. \OCTAVE is also not capable of running \EXPODE, since it currently does not support nested functions which are used regularly in the toolbox. At\\[3pt]
\url{http://www.am.uni-duesseldorf.de/en/Research/03_Software.php}\\[3pt]
you can download the packages containing the \EXPODE toolbox. Two different versions are available, a package for users and an extended one for developers. The latter contains some additional tools helpful for extending \EXPODE. Usually the user package should be sufficient.

To install, just unpack the archive obtained from the above link. In a UNIX environment you can do this by typing
\begin{sniplet}{\$}
	tar xzvf expode-VERSION.tgz
\end{sniplet}
in the directory where you downloaded the package. This will give you an \path{expode} subdirectory. To make it available in \MATLAB, just add the package's root to \MATLAB's path and run the \code{initPaths} function with
\begin{sniplet}{>{>}}
	addpath \path{/download/path/expode}; \\
	initPaths;
\end{sniplet}
To make a permanent installation for the current user, put the above line into your \file{startup.m} file. See \MATLAB's help for more information.
% 
% \begin{itemize}
% 	\item need a recent version of \MATLAB ($\geq$ 7.2), ask for feedback if runs on older versions ($\geq$ 7.0), cannot run on 6.5 due to lack of function handles, no Octave support
% 	\item the \EXPODE package is available at ....
% 	\item unzip
% 	\item add to \MATLAB path
% \end{itemize}

\subsection{Quick Start}

To get a first impression of \EXPODE we start with running some of the included examples. %Suppose we are in the \MATLAB shell at the \path{/download/path/expode} path. Run the following lines
To access the examples we add the examples directory to the \MATLAB path. Run the following lines
\begin{sniplet}{>{>}}
% 	cd \path{examples}; \\
	addpath \path{/download/path/expode/examples}; \\
	[\var{t}, \var{y}] = Heat1D([], [], 'run');
\end{sniplet}
to solve a Heat equation with a time-dependent source term in one dimension. Use
\begin{sniplet}{>{>}}
	help Heat1D;
\end{sniplet}
to get some more information on the equation. The solution over time will be plotted in a mesh plot. All examples contained in the package can be run by simply calling them without arguments. Short information on the equations are contained in their helptexts. To experiment with solver it is convenient that we now run the example manually.
\begin{sniplet}{>{>}}
% 	\c{switch directory from \path{/download/path/expode}} \\ %, add paths} \\
% 	cd \path{examples}; \\
% 	%addpath('..'); \\
% 	\\
	\c{paramerters} \\
	\var{epsilon} = 0.1; \var{gamma} = 0.1; \var{N} = 100; \\
	\\
	\c{get initial conditions} \\
	[\var{tspan}, \var{y0}, \var{options}] = Heat1D([], [], \str{init}, \var{epsilon}, \var{gamma}, \var{N}); \\
	\\
	\c{run the example} \\
	[\var{t}, \var{y}] = \EXPODECMD(@Heat1D, \var{tspan}, \var{y0}, \var{options});
\end{sniplet}
Now we can start playing with options and parameters. Switching to direct the solver for the matrix functions, we use
\begin{sniplet}{>{>}}
	\var{options} = \EXPSET(\var{options}, \str{MatrixFunctions}, \str{direct});
\end{sniplet}
Other options are set similarly. To get interactive help for options, we use one of the integrator specific set commands, i.e.
\begin{sniplet}{>{>}}
	\var{options} = \EXPRBSET(\var{options}, \str{MatrixFunctions}, \str{direct}); \\
	\var{options} = \EXPRBSET(\var{options}, \str{MinStep}, -1);
\end{sniplet}
This will do some checks on the values set for an option. The second line will throw an error, because step sizes must be positive. An overview of the available options for an integrator is available by calling the integrator info without arguments. More detailed information on a specific option can be shown with this command as well:
\begin{sniplet}{>{>}}
	\EXPRBINFO~~~~~~~~~\c{prints all available options for exprb} \\
	\EXPRBINFO MinStep \c{prints helptext for MinStep option}
\end{sniplet}
The following code will produce a logarithmic error vs. time step plot. We need to solve an equation where we can evaluate an exact solution. This is done by calling \code{ode(t, [], 'exact')}. We use the \code{semi1} example here, see its helptext for information on the equation. \code{'exprb'} indicates to only plot results for the Rosenbrock-type methods, the \code{''} uses the direct solver for the matrix functions and $\code{N} = 50$ is the number of gridpoints for the spatial discretization of \code{semi1}'s underlying partial differential equation.
\begin{sniplet}{>{>}}
	allMethods(@semi1, \str{exprb}, \str{{}}, [], 50);
\end{sniplet}
The easiest way to implement an own differential equation is to start with some tutorial files in \path{examples/Hello_World}. Then you can modify example files in the \path{examples} directory.
% The easiest way to implement an own differential equation to solve is by modifying one of the example files in the \path{examples} directory. 
\code{MinEx.m} or \code{Template.m} contain a lot of comments for an easy start. \code{MinEx.m} is simpler, where as \code{Template.m} uses more advanced features.
% The following code will produce a logarithmic error vs. step size plot. We need to solve an equation where we know an exact solution. Therefore we chose the \code{semi1} example, see its helptext for information on the equation. To save memory we use a dense output generator. Since we only need the numerical solution at the end, which is returned exactly by a Hermite interpolation formula, we do not loose accuracy.
% \begin{sniplet}{>{>}}
% 	\c{We need to modify the options, we set N=50 gridpoints in space} \\
% 	[ tpsan, y0, opts ] = semi1([], [], 'init', 50); \\
% 	\c{We don't want plots} \\
% 	opts = exprbset(opts, 'OutputFcn', 'off'); \\
% 	\c{Set dense output generator} \\
% 	opts = exprbset(opts, 'DOGenerator', 'hermite'); \\
% 	\c{Set matrix function evaluator to direct method} \\
% 	opts = exprbset(opts, 'MatrixFunctions', 'direct'); \\
% 	\c{For errorPlot: N = 50} \\
% 	opts.VARARGIN = \{ 50 \}; \\
% 	\c{Run the error plot}
% 	\errorPlot(@semi1, \{ opts \});
% \end{sniplet}
% \begin{itemize}
% 	\item cd examples
% 	\item execute examples
% 	\item set options
% 	\item errorPlot, allMethods, ... commands
% \end{itemize}

%% file: docutext.tex
\ifdefined\mathe
		\Abschnitt{The \code{exprb} Integrator}{exprb}
			The \EXPRB integrator solves a given ordinary differential equation with an exponential Rosenbrock-type method \eqref{Rosenbrock-scheme}.
			\def\INTCMD{\EXPRB}
			\def\SETCMD{\EXPRBSET}
	\else
		\Abschnitt{The \code{expode} Integrator}{expode}
			The \EXPODE integrator package solves a given ordinary differential equation of the form
			\[
				y' = F(t, y) = A y + g(t, y)
			\]
			with a number of exponential time integration methods, where the second form of the equation is only required for the semilinear solvers. For details on a specific method, see one of the later sections. The option \refOption{integrator} is used to determine the actual integrator to use, see sections \refAbschnitt{expset} and \refAbschnitt{expodeopts} for details on how to set und use it.
			\def\INTCMD{\EXPODECMD}
			\def\SETCMD{\EXPRBSET}
			
	\fi
	After downloading and extracting the code package, you will get a directory called \file{expode}.
	\ifdefined\mathe
		This directory contains the \EXPRB integrator, the options helper \code{exprbset}, the information command \code{exprbinfo},
	\else
		This directory contains the \EXPODECMD integrator and the options helper \code{expset}. %, the information command \code{exprbinfo}
		This is the general integrator, which calls the appropriate integrator as needed. Then for each supported integrator class there exists a dedicated integrator command (e.g. \EXPRB), an options helper that understands the specific options for this integrator (e.g. \code{exprbset}) and an info command (e.g. \code{exprbinfo}) to support the user. Furthermore, there is
	\fi
	an analogon to \MATLAB's \DEVAL, \DEVALEXP and the function \code{initPaths} to setup \MATLAB's searchpath. Additionally there are some subdirectories. One of these is the \file{example} directory. The programs in that directory can be run directly and be used as a reference how to call the integrator.
	
	To use \INTCMD you either have to move the package's content to your \MATLAB working directory or use
%To use \EXPRB, you either have to move all files and directories from the \code{expkit} directory into you working directory or use
%, a helper function \code{initIntegrator.m} to setup paths needed to use \EXPRB and subdirectories. To use \EXPRB, you either have to move all files and directories from the \code{expode} directory into you working directory or use
	\begin{Matlab}{addpath}%
		addpath /path/to/expode
	\end{Matlab}
	to tell \MATLAB where to find the integrator. You can also use
	\begin{Matlab}{addexamplepath}%
		addpath /path/to/expode/examples
	\end{Matlab}
	to directly run one of the examples.
	
	\INTCMD's user interface is adapted to \MATLAB's internal integrators. A call with all available arguments is of the following form:
	\begin{Matlab}{\ifdefined\mathe exprbfull\else expodefull\fi}%
		[\v{t}, \v{y}] = \INTCMD{}(\v{@\ODE{}}, \opt{\v{@jac}}, \opt{\v{tspan}, \v{y_0}, \v{opts}}, \opt{\v{varargin}});
	\end{Matlab}
	To keep the explanation well-arranged, the possible call combinations will be explained step by step. The simplest way to invoke the integrator is
	\begin{Matlab}{exprb}%
		[\v{t}, \v{y}] = \INTCMD{}(\v{@\ODE{}}, \v{tspan}, \v{y_0});
	\end{Matlab}
	where \ODE is a function to evaluate the differential equation (and its first derivative) at given data \var{y} and \var{t}. \code{\var{tspan} = [\var{t_0}, \var{T}]} defines the integration interval. You can also use \code{\var{tspan} = [\var{t_0}, \dots, \var{t_N}]}
	\ifx\mathe\undefined
		for some integrators
	\fi
	instead. This will result in a solution evaluated at the given times instead of the ones chosen by the integrator internally. In that case, clearly we have \code{\var{t} = \var{tspan}}. This kind of solution is called \emph{dense output} because one would typically use \code{\var{tspan} = \command{LINSPACE}(\var{t_0}, \var{T}, N + 1)} with different output data than the step selection would produce.
	\ifx\mathe\undefined
		\EXP4 allows dense output directly as it has an own dense output formula. For all other integrators there is only a generic hermite interpolation formula which is not suitable for stiff problems. This is disabled by default, but can be switched on via the \refOption{DOGenerator} option. Use this feature with care.
	\fi
	\var{y_0} is the initial condition to solve the equation with.
	
	The \ODE function has to be callable in the following way:
	\begin{Matlab}{dglfun}%
		\v{res} = \ODE{}(\var{t}, \var{y}, \opt{\var{flag}});
	\end{Matlab}
	\var{t} and \var{y} are the time and phase space variables respectively. 
	\ifdefined\mathe
		\var{flag} can either be omitted or given as \code{\var{flag} = \flag{jacobian}} or \code{\var{flag} = \flag{df\_dt}}. In the first case \ODE should return the evaluation of the right hand side of the equation at \var{t} and \var{y}, in the second
	\else
		\var{flag} controls the output of the function. An omitted or empty \var{flag} should trigger the evaluation of the right hand side. If the integrator used is for semilinear problems (\EXPRK and \EXPMSSEMI) \code{\var{flag} = \flag{linop}} should trigger the evaluation of the linear part (the matrix). In case of the linearized integrators \EXPRB, \EXPMS and \EXP4, \var{flag} can be given as \code{\var{flag} = \flag{jacobian}} or \code{\var{flag} = \flag{df\_dt}}. In the first
	\fi
	case the derivative of the right hand side with respect to \var{y} -- here called Jacobian, see the introduction -- has to be evaluated
	\ifdefined\mathe
		for the same data.
	\else
		 at \var{t} and \var{y}.
	\fi
	The last case is only needed if the differential equation is non-autonomous. Then \ODE should return the derivative of the right hand side with respect to the time variable \var{t}. If the differential equation is non-autonomous, you have to set the option \refOption{NonAutonomous} to \optVal{on}, see the options section. A typical function body for
	\ODE is presented
	\ifx\mathe\undefined
		for a linearized integrator 
	\fi
	below.
	\begin{Matlabfun}{dglfunsample}%
		\FUNCTION \v{res} = \ODE{}(\v{t}, \v{y}, \v{flag}) \\%
			\IF \v{nargin} == \l{2} || \ISEMPTY{}(\v{flag}) \\%
				\var{flag} = \s{}; \\%
			\END \\%
			\\%
			\SWITCH \v{flag}\\%
				\CASE \str{{}}\\%
					\v{res} = evaluation of the right hand side; \\%
				\CASE \str{jacobian} \\%
					\v{res} = evaluation of the Jacobian; \\%
				\CASE \str{linop} \\%
					\v{res} = evaluation of the linear part; \\%
				\CASE \str{df\_dt} \\%
					\v{res} = evaluation of the derivative of the ... \\%
					~~~~~~right hand side w.r.t.~\v{t}; \\%
				\OTHERWISE \\%
					\ERROR{}(\str{Unknown flag:~\%s.}, flag); \\%
			\ENDCASE \\%
			\\%
			\RETURN \v{res}; \\%
		\END
	\end{Matlabfun}
	\ifx\mathe\undefined
		If you want to use both integrator types, simply add cases for \str{jacobian}, \str{df\_dt} \emph{and} \str{linop}.
	\fi
	The \code{otherwise} case ist very helpful to find errors. All function handles used can be replaced by inlines or the functions' names as strings. All of these can be evaluated by the \EXPODE integrators.
	
	It is possible to source out the evaluation of the Jacobian
	\ifx\mathe\undefined
		or the linear part
	\fi
	into an additional function. In that case, \INTCMD has to be called this way:
	\begin{Matlab}{exprbjac}%
		[\v{t}, \v{y}] = \INTCMD{}(\v{@\ODE{}}, \v{@jac}, \v{tspan}, \v{y_0});
	\end{Matlab}
	where \code{\command{function} \var{J} = \var{jac}(\var{t}, \var{y})} is the new function for the Jacobian.
	\ifx\mathe\undefined
		The same applies to an evaluation function of the linear part, which will be called \var{jac} as well for simplicity. Note that for consistency reasons this function should allow parameters \var{t} and \var{y} as well, which can be ignored by the code.
	\fi
	The handle to
	\ifdefined\mathe
		the Jacobian
	\else
		these functions
	\fi
	can alternatively be given as an integrator option. See the \refOption{Jacobian} option in section 
	\ifdefined\mathe
		\refAbschnitt{exprbopts}
	\else
		\refUnterabschnitt{linearizationopts} or the \refOption{LinOp} option in section \refUnterabschnitt{semilinopts}
	\fi
	for details. If
	\ifdefined\mathe
		this option is
	\else
		these options are
	\fi
	switched \str{off} and \refOption{JacobianV}
	\ifdefined\mathe
		is
	\else
		or \refOption{LinOpV} are
	\fi
	\str{on} instead, \var{jac} will be interpreted as a function to evaluate the product of the Jacobian
	\ifx\mathe\undefined
		or linear part
	\fi
	with a vector. In that case, \code{\command{function} \var{res} = \var{jac}(\var{t}, \var{y}, \var{v})} has to be callable.
	
	To use integrator options, invoke \INTCMD in the following way:
	\begin{Matlab}{exprbopts}%
		[\v{t}, \v{y}] = \INTCMD{}(\v{@\ODE{}}, \opt{\v{@jac}}, \v{tspan}, \v{y_0}, \var{opts});
	\end{Matlab}
	The construction of such an options object will be discussed in the \hrefAbschnitt{expset}{next section}.
	
	If you want to hand over additional parameters to the \ODE function -- stiffness parameters for springs or some dimensions for example -- these parameters can be passed to the integrator as follows:
	\begin{Matlab}{exprbparam}%
		[\v{t}, \v{y}] = \INTCMD{}(\v{@\ODE{}}, \opt{\v{@jac}}, \v{tspan}, \v{y_0}, \v{opts}, \var{varargin});
	\end{Matlab}
	\code{varargin} can be an arbitrary number of additional argument separated by comma. These parameters will be passed to the \ODE function directly via
	\begin{Matlab}{dglfunparam}%
		\v{res} = \ODE{}(\var{t}, \var{y}, \var{flag}, \var{varargin}\opt{:});
	\end{Matlab}
	Please \emph{do not} use a \command{struct} object as the first of the \var{varargin} arguments if you do not use integrator options. The integrator would confuse this argument with an options structure.
	
% 	Also please \emph{do not} use a \command{char} array (string) as the first of the \var{varargin} arguments, because the \ODE function would confuse this with the flag. A skeleton for the \ODE function with additional arguments could look like
	
	If you still want to do so, you have to provide an empty vector (\code{[]}) as the options parameter.
	
	\begin{Matlabfun}{dglfunparamssample}%
		\FUNCTION \v{res} = \ODE{}(\v{t}, \v{y}, \v{flag}, \v{varargin}) \\%
			\IF \v{nargin} == \l{2} || \ISEMPTY{}(\v{flag}) || \~\ISCHAR{}(\v{flag}) \\%
				\IF \~\ISCHAR{}(\v{flag}) \\%
					\v{varargin} = \{\v{flag}, \v{varargin}\{:\}\}; \\%
				\END \\%
				\var{flag} = \s{}; \\%
			\END \\%
			\\%
			\c{The rest can be the same as \refSnipplet{dglfunsample}{before}.} \\%
		\END
	\end{Matlabfun}
	
	If \var{varargin} is used, it will also be passed to the \var{jac} function if it is available. It will be called with \code{J = \var{jac}(\var{t}, \var{y}, \var{varargin}\{:\})} in that case.
	
	It is also possible to omit the integration interval, the initial condition, and the integrator options if the \ODE function is able to provide them. If \var{varargin} is used, you need to use an empty vector (\code{[]}) for the omitted parameters. Calling \INTCMD with
	\begin{Matlab}{exprbdglinits}%
		[\v{t}, \v{y}] = \INTCMD{}(\v{@\ODE{}}, \opt{\v{@jac}});
	\end{Matlab}
	or
	\begin{Matlab}{exprbdglinitsvarargin}%
		[\v{t}, \v{y}] = \INTCMD{}(\v{@\ODE{}}, \opt{\v{@jac}}, [], [], [], \opt{\v{varargin}});
	\end{Matlab}
	will invoke \ODE with the \str{init} flag:
	\begin{Matlab}{dglinits}%
		[\v{tspan}, \v{y_0}, \var{opts}] = \ODE{}([], [], \s{init}, \opt{\v{varargin}});
	\end{Matlab}
	The example code for the \ODE function has to be extended by another \command{case} for \str{init}, as can be seen here.
	\begin{Matlabfun}{dglinitscode}%
		\ind\ind\ind\indent\CASE \str{init} \\%
			\v{res}\{1\} = [\var{t_0}, \var{T}]; \\%
			\v{res}\{2\} = \var{y_0}; \\%
			\v{opts} ~~= \SETCMD{}(\v{option1}, \v{value1}, \v{option2}, \v{value2}, ...); \\%
			{}~~~~~~~\c{see section \refAbschnitt{expset}, \refSnipplet{exprbsetsample}{exprbset}} \\%
			\v{res}\{3\} = \var{opts};
	\end{Matlabfun}
	
	To be able to evaluate the solution at an arbitrary time after the integration is finished, the call
	\begin{Matlab}{dglsolveroutput}%
		\var{sol} = \EXPODE{}(\v{@\ODE{}}, \opt{\v{@jac}}, \opt{\v{tspan}, \v{y_0}, \v{opts}}, \opt{\v{varargin}});
	\end{Matlab}
	will generate a variable \var{sol} which can be used with the \DEVALEXP function via
	\begin{Matlab}{devalexp}%
		[ \var{y}, \var{dydt} ] = \DEVALEXP{}(\var{sol}, \var{t});
	\end{Matlab}
	The argument \var{t} then is a vector of times in the integration interval. Then \var{y} contains the evaluation of the numerical solution and \var{dydt} its first derivative ot the times in \var{t}. Depending on the dense output generator used, \var{sol} can get quite large. For the \EXP4 integrator it has to save all inner stages for each integration step, that is a total eight vectors per step. Note that this only works if a dense output formula is available.
	
\Abschnitt{Using Options}{expset}
	\ifdefined\mathe
		\SETCMD is used to set options for the \INTCMD integrator.
	\else
		For each of our integrators we have a command to set options. Exemplary we show their usage with \SETCMD, the one for the \EXPRB integrator. 
	\fi
	It is an extension to \MATLAB's \ODESET. It would have been preferable to use \ODESET directly, but it has a fixed set of available options and is not extensible to support a larger set.
	
	Both \ODESET and \SETCMD create a so called \emph{options structure} or \emph{options object}, which contains the options set by the user. These structures can be passed to the integrators. It has been explained in the \refSnipplet{\ifdefined\mathe exprb\else expode\fi}{previous section} how to do that for \INTCMD.
	
	\SETCMD is compatible with the option objects created by \ODESET, so one can create an options structure with \ODESET and then extend it with \SETCMD. Unfortunately, it is not possible the other way around, since \ODESET removes all options it does not know about. The first sequence is more important though. Usually, a program will be written using \MATLAB's standard tools and then be extended to use other integrators.
	
	The \SETCMD function can be called analogously to \ODESET. The simplest way to use it is
	\begin{Matlab}{expsetsample}%
		\v{opts} = \SETCMD{}(\v{option}, \v{value});
	\end{Matlab}
	where \var{option} is an option's name and \var{value} the user's choice. The latter one can be scalar, vector valued, logical, a string, or a function, depending on the option chosen. All supported option types are listed later in this section.
	\ifdefined\mathe
		The \refSnipplet{exprbinfo}{\EXPRBINFO{}}% command provided by the \EXPODE package can print options supported by the \EXPRB integrator.
	\else
		The info commands such as \EXPRBINFO{}
	\fi
	provided by the \EXPODE package can print options supported by the \ifdefined\mathe\EXPRB\else\EXPODE \fi integrator\ifx\mathe\undefined s\fi.
	
	It is also possible to pass more option-value-pairs at once and extend already existing option objects:
	\begin{Matlab}{expsetsamplemulti}%
		\v{opts} = \SETCMD{}(\v{option1}, \v{value1}, \v{option2}, \v{value2}, ...); \\%
		\v{opts} = \SETCMD{}(\v{opts}, \v{option3}, \v{value3}, ...);
	\end{Matlab}
	
	\ifdefined\mathe
		For a more comfortable usage, \EXPRBSET checks if the provided options are valid.
	\else
		For a more comfortable usage, you can use the integrator specific set commands (\EXPRBSET, \EXPRKSET, etc.) to check if the provided options are valid. \EXPSET itself cannot do this, since it does not know which integrator will be used. You can tell it by using
		\begin{Matlab}{expsetsampleintname}%
			\v{opts} = \SETCMD{}(\v{intname}, \opt{\v{opts}}, \v{option1}, \v{value1}, \v{option2}, \v{value2}, ...);
		\end{Matlab}
		where \var{intname} is the name of the integrator, e.g. \str{exprb}, \str{exprk}.
	\fi
	
	In the remaining part of this section, the supported value types for options will be discussed. All of the listed types can be checked for validity by \ifdefined\mathe\EXPRBSET\else the set commands\fi. Which type is allowed for which \ifdefined\mathe\EXPRB\fi option is described in the \hrefAbschnitt{exprbopts}{next sections}. The \ifx\mathe\undefined info commands \else\EXPRBINFO command \fi in \MATLAB can provide this information.
	\begin{Matlab}{exprbinfo}%
		\EXPRBINFO{}(\opt{optname});
	\end{Matlab}
	will display a list of all options supported by \ifdefined\mathe the integrator\else\EXPRB\fi if \var{optname} is omitted. If it is given, more detailed help on the specific option will be printed. Using \code{\var{optname} = \str{-}} will give this detailed help on all options. Warning: the output is quite long, it should be used with
	\begin{Matlab}{moreexprbinfo}%
		\MORE \ON; \EXPRBINFO{}(\s{-}); \MORE \OFF;
	\end{Matlab}
	
	The output of \ifdefined\mathe\EXPRBINFO\else the info command \fi contains the option's type in the first line in squared brackets. For the \refOption{AbsTol} option this line reads
	
	\code{AbsTol\! -\! Absolute\! error\! tolerance\! [\! positive\! scalar\! |\! positive\! vector\! \{1e-06\}\! ]}.
	
	The vertical bar separates alternative types, so the type is \refOptType{positive} \refOptType{scalar} or \refOptType{positive} \refOptType{vector}. The default value is set in braces after the type. The two types are each a combination of two other types: \refOptType{positive} and \refOptType{scalar} or \refOptType{vector} respectively. Which types are compatible will be explained below. Another example for types is given by the \refOption{JacobianV} option. Here the \EXPRBINFO command prints
	
	\code{JacobianV - ... [ function\_handle | \{'off'\} | 'on' ]}.
	
	The allowed types in this case are \refOptType{boolean} or \refOptType{function\_handle}. For types \refOptType{boolean} and its generalization \refOptType{list}, a list of accepted values will be stated. For a \refOptType{boolean} this list is fixed and consists of the values \str{on} (\code{true}) and \str{off} (\code{off}), for the \refOptType{list} type the list depends on the option. The listed values are always quoted strings. Now all option types will be discussed in detail.
	\begin{OptionType}{scalar}
		This type can be any numerical type that is supported by \MATLAB. \opttype{scalar} can be combined with \refOptType{integer}, \refOptType{positive}, \refOptType{non-negative}, \refOptType{negative} and \refOptType{non-positive}.
	\end{OptionType}
	\begin{OptionType}{vector}
		This type can be a vector of any numerical type that is supported by \MATLAB. \refOptType{scalar} is the special case of this type with length one, so any scalar is a vector as well. \opttype{Vector} can be combined with \refOptType{integer}, \refOptType{positive}, \refOptType{non-negative}, \refOptType{negative} and \refOptType{non-positive}.
	\end{OptionType}
	\begin{OptionType}{matrix}
		This type can be a matrix of any numerical type that is supported by \MATLAB. \refOptType{scalar} and \refOptType{vector} are special cases of this type, so any scalar and vector is a matrix as well. \opttype{Matrix} can be combined with \refOptType{integer}, \refOptType{positive}, \refOptType{non-negative}, \refOptType{negative} and \refOptType{non-positive}.
	\end{OptionType}
	\begin{OptionType}{integer}
		This type can be any integer number. This does not enforce the use of \MATLAB's \literal{int8}, \literal{int16} or similar types. \var{value} is accepted as an \opttype{integer} as long as the expression \code{(\var{value} - \command{round}(\var{value})) == 0} ist \code{true}. \opttype{integer} can be combined with \refOptType{vector}, \refOptType{scalar}, \refOptType{positive}, \refOptType{non-negative}, \refOptType{negative} and \refOptType{non-positive}.
	\end{OptionType}
	\begin{OptionType}{positive}
		This type can be any numerical type that is supported by \MATLAB (scalar, vector or matrix) where all components are strictly positive. \opttype{positive} can be combined with \refOptType{matrix}, \refOptType{vector}, \refOptType{scalar} and \refOptType{integer}.
	\end{OptionType}
	\begin{OptionType}{non-negative}
		This type allows any numerical \MATLAB type (scalar, vector or matrix) where all components are not negative (positive or zero). \opttype{non-negative} can be combined with \refOptType{matrix}, \refOptType{vector}, \refOptType{scalar} and \refOptType{integer}.
	\end{OptionType}
	\begin{OptionType}{negative}
		This type can be any numerical type that is supported by \MATLAB (scalar, vector or matrix) where all components are strictly negative. \opttype{negative} can be combined with \refOptType{matrix}, \refOptType{vector}, \refOptType{scalar} and \refOptType{integer}.
	\end{OptionType}
	\begin{OptionType}{non-positive}
		This type allows any numerical \MATLAB type (scalar, vector or matrix) where all components are not positive (negative or zero). \opttype{non-positive} can be combined with \refOptType{matrix}, \refOptType{vector}, \refOptType{scalar} and \refOptType{integer}.
	\end{OptionType}
	\begin{OptionType}{index}
		This type is an abbreviation for \refOptType{positive} \refOptType{integer} \refOptType{scalar}. \opttype{index} can therefore not be combined with any other type.
	\end{OptionType}
	\begin{OptionType}{indices}
		This type is an abbreviation for \refOptType{positive} \refOptType{integer} \refOptType{vector}. \opttype{indices} therefore cannot be combined with any other type.
	\end{OptionType}
	\begin{OptionType}{boolean}
		The values of this type can be given in different forms. Either use \str{on} and \str{off} or use the numerical values \literal{true} and \literal{false}. \opttype{boolean} is not combinable with any other type. Alternative string values are \str{yes} and \str{no} or \str{true} and \str{false}.
	\end{OptionType}
	\begin{OptionType}{list}
		The values for this type can be given in different forms. Either use one of the available values printed by the \ifdefined\mathe\EXPRBINFO\else info \fi command or use the number of the list entry, start counting from \literal{0}. Example: The \refOption{MatrixFunctions} option's info states (ignore \code{function\_handle} here, it will be discussed below)
		
		\code{MatrixFunctions - ... [ \{\str{direct}\} | \str{arnoldi} | function\_handle ]}
		
		which permits to set the values \str{direct} and \str{arnoldi}, in that order. Therefore \str{direct} has number \literal{0} and \str{arnoldi} has number \literal{1}. Some special \opttype{list} type options define specific numeric values, which will be stated in parentheses behind the string value. See, for instance, option \refOptionPrefix{exprb}{Order} of the \EXPRB integrator. It is possible to set the value for the option \refOption{MatrixFunctions} to use Arnoldi's method in the following ways:
		
		\begin{Matlab}{expsetlist}%
			\v{opts} = \SETCMD{}(\s{MatrixFunctions}, \s{arnoldi}); \c{or equivalently} \\%
			\v{opts} = \SETCMD{}(\s{MatrixFunctions}, \l{1});
		\end{Matlab}
		
		\opttype{list} cannot be combined with any other type.
	\end{OptionType}
	\begin{OptionType}{text}
		This type can be any \MATLAB \command{char} array (string). \opttype{text} cannot be combined with any other type.
	\end{OptionType}   
	\begin{OptionType}{struct}
		This type can be any \MATLAB \command{struct} object. \opttype{struct} cannot be combined with any other type.
	\end{OptionType}
	\begin{OptionType}{function\_handle}
		This type can be either a \MATLAB function handle, a \MATLAB inline function or a \command{char} array (string). A function handle to a function \code{func} is generated by \code{@func}. An inline function is generated by \MATLAB's \command{inline} command. If the argument is a string, it has to represent a function's name: use \str{func} for the function \code{func}. \opttype{function\_handle} cannot be combined with any other type.
	\end{OptionType}

\ifdefined\mathe
	\input{exprboptions}
\else
	\input{expodeoptions}

\fi

\input{writingMatrixFunctions}

	\input{writingexpode}

%% file: expodeoptions.tex
\Abschnitt{Available Options}{expodeopts}
	
	In the following sections, the options available for the \EXPODE integrators will be discussed. The first option discussed will be an exception to this, as it is used to choose the actual integrator when calling expode directly. The rest of this section will be split into several subsections, starting with options common to all integrators, followed by options for semilinear, then linearized integrators. The last two subsections are dedicated to constant and variable step size integrators. Let's start with the integrator option.

\begin{Option}{}{Integrator}%
	\input{options/Integrator.tex}\end{Option}%

	\Unterabschnitt{Options Common to All Integrators}{expodeopts}
		
		Here we describe options common to all integrators. The options will be sorted into several groups. First some options for controlling the integration process, then some properties for the ODE followed by output control and debugging options.

\begin{Option}{}{AbsTol}%
	\input{options/AbsTol.tex}\end{Option}%

\begin{Option}{}{RelTol}%
	\input{options/RelTol.tex}\end{Option}%

\begin{Option}{}{NormControl}%
	\input{options/NormControl.tex}\end{Option}%

\begin{Option}{}{MatrixFunctions}%
	\input{options/MatrixFunctions.tex}\end{Option}%

\begin{Option}{}{KrylovTestIndex}%
	\input{options/KrylovTestIndex.tex}\end{Option}%

		Properties of the differential equation can be set by the following options.

\begin{Option}{}{NonAutonomous}%
	\input{options/NonAutonomous.tex}\end{Option}%

\begin{Option}{}{Complex}%
	\input{options/Complex.tex}\end{Option}%

\begin{Option}{}{Structure}%
	\input{options/Structure.tex}\end{Option}%

\begin{Option}{}{GFcn}%
	\input{options/GFcn.tex}\end{Option}%

		Now we will discuss how to control the output of the integrator.

\begin{Option}{}{DOGenerator}%
	\input{options/DOGenerator.tex}\end{Option}%

\begin{Option}{}{Refine}%
	\input{options/Refine.tex}\end{Option}%

\begin{Option}{}{OutputFcn}%
	\input{options/OutputFcn.tex}\end{Option}%

\begin{Option}{}{OutputSel}%
	\input{options/OutputSel.tex}\end{Option}%

		The last set of options in this section deals with debugging and logging.

\begin{Option}{}{Stats}%
	\input{options/Stats.tex}\end{Option}%

\begin{Option}{}{MatrixFunctionStats}%
	\input{options/MatrixFunctionStats.tex}\end{Option}%

\begin{Option}{}{Waitbar}%
	\input{options/Waitbar.tex}\end{Option}%

\begin{Option}{}{ClearInternalData}%
	\input{options/ClearInternalData.tex}\end{Option}%

	\Unterabschnitt{Options for Semilinear Integrators}{semilinopts}
		
		In this subsection we will discuss options common to the semilinear integrators \EXPRK and \EXPMSSEMI. To evaluate the linear operator we have the following two options.

\begin{Option}{}{LinOp}%
	\input{options/LinOp.tex}\end{Option}%

\begin{Option}{}{LinOpV}%
	\input{options/LinOpV.tex}\end{Option}%

		Additionally you can control the logging output with this option.

\begin{Option}{}{LinOpStats}%
	\input{options/LinOpStats.tex}\end{Option}%

	\Unterabschnitt{Options for Linearized Integrators}{linearizationopts}
	
		In this subsection we will discuss options common to the linearized integrators \EXPRB, \EXPMS and \EXP4. Again the options are split into smaller groups to keep the information clearly presented. The first couple of options deal with the evaluation of the needed data. Additionally we have one additional option for the equation's properties and one for logging.

\begin{Option}{}{Jacobian}%
	\input{options/Jacobian.tex}\end{Option}%

\begin{Option}{}{JacobianV}%
	\input{options/JacobianV.tex}\end{Option}%

\begin{Option}{}{GJacobian}%
	\input{options/GJacobian.tex}\end{Option}%

\begin{Option}{}{GJacobianV}%
	\input{options/GJacobianV.tex}\end{Option}%

		The differential equation can now have this additional property.

\begin{Option}{}{Semilin}%
	\input{options/Semilin.tex}\end{Option}%

		Please note that for the semilinear integrators in the next section this option is always assumed \str{on}.
		
		The last option extends the logging capabilities.

\begin{Option}{}{JacobianStats}%
	\input{options/JacobianStats.tex}\end{Option}%

	\Unterabschnitt{Options for Constant Step Size Integrators}{conststepopts}
		
		Now we can move on to integrators without step size control, that currently are \EXPRK, \EXPMS and \EXPMSSEMI. Since there is not so much to control, we only have one option.

\begin{Option}{}{StepSize}%
	\input{options/StepSize.tex}\end{Option}%

	\Unterabschnitt{Options for Variable Step Size Integrators}{varstepopts}
		
		The last general set of options is for the variable step size integrators \EXPRB and \EXP4. Here we have the following options.

\begin{Option}{}{hConstant}%
	\input{options/hConstant.tex}\end{Option}%

\begin{Option}{}{InitialStep}%
	\input{options/InitialStep.tex}\end{Option}%

\begin{Option}{}{MaxStep}%
	\input{options/MaxStep.tex}\end{Option}%

\begin{Option}{}{MinStep}%
	\input{options/MinStep.tex}\end{Option}%

		And we have another option for logging again.

\begin{Option}{}{StepStats}%
	\input{options/StepStats.tex}\end{Option}%

\Abschnitt{\code{exprk} -- The Exponential Runge-Kutta Integrator}{exprkopts}
	
	Since we are now familiar with the common options, it is now time to move to the actual integrators. This section will start with the exponential Runge-Kutta methods. \EXPODE only supports explicit schemes. A good overview over this class can be found in \cite{HocO05SIAM}. \EXPRK is a semilinear solver and only supports constant step size integrators. This can easily be extended to variable step size if appropriate error estimators are available for the schemes.
	
	To call \EXPRK, use the following syntax:
	\begin{Matlab}{exprkfull}%
		[\v{t}, \v{y}] = \EXPRK{}(\v{@\ODE{}}, \opt{\v{tspan}, \v{y_0}, \v{opts}}, \opt{\v{varargin}});
	\end{Matlab}
	See section \refAbschnitt{expode} for details on the arguments.
	
	The \EXPRK integrator supports a number of Runge-Kutta schemes. See \cite{HocO05SIAM} for a reference. The scheme to use can be selected via the \refOption{Scheme} option, see below. It is possible to either select one of the predefined schemes or even provide your own. To do so, you have to create a scheme with the \code{rkScheme} function from the \path{exprk} subdirectory of the \EXPODE distribution. \code{rkScheme} has quite a detailed help built-in (\code{help rkScheme} in \MATLAB). Use \code{initPaths(\str{exprk})} to put the \code{exprk} directory into your \MATLAB-Path. To better understand the usage of this tool, its usage is shown by constructing Krodstad's \cite{Kro05} scheme here.
	
	\begin{Matlabfun}{rkSchemeExample}%
		\c{Make sure tp have \path{/path/to/expode} in your \MATLAB path} \\%
		initPaths(\s{exprk}); \\%
		\c{Create an empty 4-node scheme} \\%
		sc = rkScheme(4, 4); \\%
		\c{Set the nodes} \\%
		sc.c = [ 0 1/2 1/2 1 ]; \\%
		\c{U${}_2$ = u${}_n$ + 1/2 $\varphi_1$(-c${}_2$ h A) (g(U${}_1$) - A u${}_n$))} \\%
		sc = rkScheme(sc, 2, 1, 1/2); \\%
		\c{U${}_3$ = u${}_n$ + (1/2 $\varphi_1$(-c${}_3$ h A) - $\varphi_2$(-c${}_3$ h A)) (g(U${}_1$) - A u${}_n$) + ...} \\%
		\c{~~~~~~~~~~$\varphi_2$(-c${}_3$ h A) (g(U${}_2$) - A u${}_n$)} \\%
		sc = rkScheme(sc, 3, 1, [ 1/2, -1 ]); \\%
		sc = rkScheme(sc, 3, 2, [ 0, 1 ]); \\%
		\c{Stage 4} \\%
		sc = rkScheme(sc, 4, 1, [ 1, -2 ]); \\%
		sc = rkScheme(sc, 4, 3, [ 0, 2 ]); \\%
		\c{Outer stage} \\%
		sc = rkScheme(sc, 1, 'b', [ 1, -3, 4 ]); \\%
		sc = rkScheme(sc, 2, 'b', [ 0, 2, -4 ]); \\%
		sc = rkScheme(sc, 3, 'b', [ 0, 2, -4 ]); \\%
		sc = rkScheme(sc, 4, 'b', [ 0, -1, 4 ]);
	\end{Matlabfun}
	
	For reference on more complex schemes that use $\varphi_k(c_l)$ for stage $U_{ij}$ or other functions than the $\varphi$'s, see the function \code{getRKSchemeFromOption} found in the \path{exprk} directory as well.
	
	%There are two options additionally to the common, constant step size and semilinear options, that will be explained now.
	In addition to the \hrefUnterabschnitt{expodeopts}{common}, \hrefUnterabschnitt{conststepopts}{constant step size} and \hrefUnterabschnitt{semilinopts}{semilinear} options, there are two additional options to determine the Runge-Kutte scheme to use, which will be explained now.

\begin{Option}{exprk}{Scheme}%
	\input{options/exprkScheme.tex}\end{Option}%

\begin{Option}{exprk}{Parameters}%
	\input{options/exprkParameters.tex}\end{Option}%

\Abschnitt{\code{exprb} -- The Exponential Rosenbrock-type Integrator}{exprbopts}
	
	Now we turn our attention to exponential Rosenbrock-type methods \cite{ExpRBPaper}. They belong to the class of linearized integrators and have an embedded error estimator, so have an adaptive time stepping.
	
	To call \EXPRB, use the following syntax:
	\begin{Matlab}{exprbfull}%
		[\v{t}, \v{y}] = \EXPRB{}(\v{@\ODE{}}, \opt{\v{@jac}}, \opt{\v{tspan}, \v{y_0}, \v{opts}}, \opt{\v{varargin}});
	\end{Matlab}
	See section \refAbschnitt{expode} for details on the arguments.
	
	Due to their simpler order conditions there are much less schemes. In addition to the \hrefUnterabschnitt{expodeopts}{common}, \hrefUnterabschnitt{varstepopts}{variable step size} and \hrefUnterabschnitt{linearizationopts}{linearization} options there are two additional options.

\begin{Option}{exprb}{Order}%
	\input{options/exprbOrder.tex}\end{Option}%

\begin{Option}{exprb}{ErrorEstimate}%
	\input{options/exprbErrorEstimate.tex}\end{Option}%

\Abschnitt{\code{expmssemi} -- The Exponential Multistep Integrator}{expmssemiopts}
	
	This section discusses the exponential multistep integrator \EXPMSSEMI, cf. \cite{Nor69} and \cite{ExpIntPaper}. As with exponential Runge-Kutta methods it's a semilinear solver and only supports constant step sizes. In contrast to these methods, it is much harder to generalize them to variable step sizes due to their construction.
	
	To call \EXPMSSEMI, use the following syntax:
	\begin{Matlab}{expmssemifull}%
		[\v{t}, \v{y}] = \EXPMSSEMI{}(\v{@\ODE{}}, \opt{\v{tspan}, \v{y_0}, \v{opts}}, \opt{\v{varargin}});
	\end{Matlab}
	See section \refAbschnitt{expode} for details on the arguments.
	
	In addition to the \hrefUnterabschnitt{expodeopts}{common}, \hrefUnterabschnitt{conststepopts}{constant step size} and \hrefUnterabschnitt{semilinopts}{semilinear} options there is one additional option to determine the step count.
	
% 	\incOptionPrefix{expmssemi}{kStep}
	%
\begin{Option}{expmssemi}{kStep}%
	\input{options/expmssemikStep.tex}\end{Option}%

\begin{Option}{expmssemi}{StartupSteps}%
	\input{options/expmssemiStartupSteps.tex}\end{Option}%

\Abschnitt{\code{expms} -- The Exponential Linearized Multistep Integrator}{expmsopts}
	
	This section discusses the exponential linearized multistep integrator \EXPMS. It is a linearized variant of the \EXPMSSEMI and only supports constant step sizes. We have implemented the general scheme introduced by Hochbruck and Ostermann in \cite{HocO10BIT} and a scheme proposed by Tokman in \cite{Tok06}. %We have implemented two schemes, one by Tokman \cite{Tok06} and one by Hochbruck and Ostermann \cite{ExpIntPaper}.
	
	To call \EXPMS, use the following syntax:
	\begin{Matlab}{expmsfull}%
		[\v{t}, \v{y}] = \EXPMS{}(\v{@\ODE{}}, \opt{\v{@jac}}, \opt{\v{tspan}, \v{y_0}, \v{opts}}, \opt{\v{varargin}});
	\end{Matlab}
	See section \refAbschnitt{expode} for details on the arguments.
	
	In addition to the \hrefUnterabschnitt{expodeopts}{common}, \hrefUnterabschnitt{conststepopts}{constant step size} and \hrefUnterabschnitt{linearizationopts}{linearization} options there is one additional option to determine the integration scheme.

\begin{Option}{expms}{kStep}%
	\input{options/expmskStep.tex}\end{Option}%

\begin{Option}{expms}{StartupSteps}%
	\input{options/expmsStartupSteps.tex}\end{Option}%

\Abschnitt{\code{exp4} -- Exponential Integrator of Order Four}{exp4opts}
	
	This section is dedicated to the \EXP4 integrator \cite{exp4}. It belongs to the class of linearized integrators and has adaptive time stepping.
	
	To call \EXP4, use the following syntax:
	\begin{Matlab}{exp4full}%
		[\v{t}, \v{y}] = \EXP4{}(\v{@\ODE{}}, \opt{\v{@jac}}, \opt{\v{tspan}, \v{y_0}, \v{opts}}, \opt{\v{varargin}});
	\end{Matlab}
	See section \refAbschnitt{expode} for details on the arguments.
	
	Since \EXP4 has its own dense output formula, it overrides the \refOption{DOGenerator} option.

\begin{Option}{exp4}{DOGenerator}%
	\input{options/exp4DOGenerator.tex}\end{Option}%

%% file: options/Integrator.tex
\optshortdesc{Integrator to use}\\
Permitted values: \code{\{\optVal{exprb}\}}, \code{\optVal{expms}}, \code{\optVal{expmssemi}}, \code{\optVal{exprk}} or \code{\optVal{exp4}}

Select the integrator to do the actual work. Simply type the integrator's
name in the \MATLAB prompt to get more information. Instead of using
this options it is preferred to run the appropriate integrator directly.
This option mainly exists for automation purposes in scripts.

%% file: options/AbsTol.tex
\optshortdesc{Absolute error tolerance}\\
Permitted type: \refOptType{positive} \refOptType{scalar} or \refOptType{positive} \refOptType{vector} \code{\{1e-06\} }

A scalar tolerance applies to all components of the solution vector.
Elements of a vector of tolerances apply to corresponding components of
the solution vector. \refOption{AbsTol} defaults to $10^{-6} $ in all solvers.
This also controls the stopping criterion in the Krylov process if matrix
functions are evaluated this way.

\optseealso \refOption{RelTol} and \refOption{NormControl}

%% file: options/RelTol.tex
\optshortdesc{Relative error tolerance}\\
Permitted type: \refOptType{positive} \refOptType{scalar} \code{\{0.001\} }

This scalar applies to all components of the solution vector, and
it defaults to $10^{-3} $ ($0.1\%$ accuracy) in all solvers.  The
estimated error in each integration step satisfies $\| err \| \ensuremath\leq 1$ in a
norm scaled with $\var{RelTol}\cdot \max (\abs (\var{y_{n}}(i), \abs (\var{y_{n-1}}(i))) + \var{AbsTol}(i)$ in each
component, where $\var{y_{n}}$ is the numerical solution at the current, $\var{y_{n-1}}$ the one
at the previous time step. This also controls the stopping criterion in the
Krylov process if matrix functions are evaluated this way.

\optseealso \refOption{AbsTol} and \refOption{NormControl}

%% file: options/NormControl.tex
\optshortdesc{Control only 2--norm of the error}\\
Permitted values: \code{\{\optVal{off}\}} or \code{\optVal{on}}

Set this property to \str{on} to request that the solver controls the error
in each integration step to meet $\| err \| \ensuremath\leq (1/n)\cdot \max (RelTol\cdot \| y \|,AbsTol)$.
By default the solvers use a more stringent componentwise error control.

\optseealso \refOption{AbsTol} and \refOption{RelTol}

%% file: options/MatrixFunctions.tex
\optshortdesc{Evaluation method for the matrix functions}\\
Permitted values: \code{\{\optVal{direct}\}}, \code{\optVal{arnoldi}} or \refOptType{function\_handle}

Set the method to evaluate the product of the matrix function with vectors.
There are two built--in methods. The default \str{direct} method uses diagonalization
of the Jacobian to do this. If the Jacobian is too large, this will be too expensive
computionally -- use \str{arnoldi} in this case. All matrix functions will then
be approximated in Krylov subspaces, the Arnoldi method is then used to compute
a nested orthonormal basis of that space.
Another alternative is to provide a custom function to compute these results.
See the section \refAbschnitt{matFun} how this method has to work.

\optseealso \refOption{Jacobian}, \refOption{JacobianV}, \refOption{LinOp}, \refOption{LinOpV} and \refOption{KrylovTestIndex}

%% file: options/KrylovTestIndex.tex
\optshortdesc{Dimensions of the Krylov subspaces to test the residual
}\\
Permitted type: \refOptType{index} or \refOptType{vector} \refOptType{of} \refOptType{indicies} \code{\{[ 1 2 3 4 6 8 11 15 20 27 36 46 57 70 85 100 ]\} }

Sequence of ascending integers which indicate the dimensions of the
Krylov subspaces where the residual is tested. If the KrylovMaxDim option
is set to zero, the largest number is the maximum dimension, otherwise dimensions
greater than specified there will be ignored. For KrylovMaxDim itself, the
residual will also be tested. If the Krylov process does not terminate within
this range, the time step size is reduced. This option only applies if
MatrixFunctions is set to \str{arnoldi}.

\optseealso \refOption{MatrixFunctions}, \refOption{KrylovMaxDim} and \refOption{KrylovAbort}

%% file: options/NonAutonomous.tex
\optshortdesc{Specifies whether the ODE is autonomous or not}\\
Permitted values: \code{\{\optVal{off}\}} or \code{\optVal{on}}

Set to \str{on} if the differential equation is not autonomous. This option
is only required for integrators using the Jacobian. It prints a message
at integrator startup for all integrators though.

%% file: options/Complex.tex
\optshortdesc{Solution is complex}\\
Permitted values: \code{\optVal{off}} or \code{\{\optVal{on}\}}

Set to \str{off} if the solution does not have complex components. This
will set all possibly numerically generated imaginary parts of the
solution to zero.

%% file: options/Structure.tex
\optshortdesc{Structure of the Jacobian/linear part}\\
Permitted values: \code{\{\optVal{none}\}}, \code{\optVal{normal}}, \code{\optVal{symmetric}}, \code{\optVal{skewsymmetric}} or \code{\optVal{diagonal}}

If the Jacobian (or linear part) has a special structure, this can
be exploited in the matrix function evaluators.
Available properties are none -- if there is no special structure,
diagonal -- if the matrix is diagonal, symmetric -- for symmetric or
Hermitian matrices, skewsymmetric -- for skew--symmetric or skew--Hermitian
matrices and normal -- for normal matrices, that are none of the three
previous.

\aliases \refOption{Symmetry}

%% file: options/GFcn.tex
\optshortdesc{Evaluation of the nonlinear part of the ODE, for semilinear problems
}\\
Permitted type and values: \refOptType{function\_handle}, \code{\optVal{off}} or \code{\{\optVal{on}\}}

This option will only be used, if \refOption{Semilin} is \str{on}.
Set this to \str{on} if the ODE file can evaluate the nonlinear part
when used with flag \str{gfun}. Set to a function handle if an external
function can evaluate it. Setting to \str{off} will result in reverting
to the standard solving method, ignoring the special structure of the ODE.

\optseealso \refOption{Semilin}

%% file: options/DOGenerator.tex
\optshortdesc{Dense Output generator}\\
Permitted type and values: \refOptType{function\_handle}, \code{\{\optVal{default}\}} or \code{\optVal{hermite}}

Set the method to generate dense output. \str{default} is the integrator's
default method, which is none for the most integrators. \str{hermite} uses hermite
interpolation and is usable with all integrators, but is not applicable for
stiff problems! Use it for testing only.
Integrators having a designated dense output formula will override and extend
this option.

To obtain solutions at specific times \var{t0},\var{t1},...,\var{tfinal} (increasing or decreasing)
use \var{tspan} = [\var{t0} \var{t1} ... \var{tfinal}] when calling
    \begin{Matlab}{OptionDOGenerator0}%
	[\var{tout},\var{yout}]=\code{expode}(\ODE{},\var{tspan},\var{y0},...);
\end{Matlab}
or
    \begin{Matlab}{OptionDOGenerator1}%
	[\var{tout},\var{yout}]=\code{expode}(\ODE{},\var{jac},\var{tspan},\var{y0},...);\\%
    \var{sol}=\code{expode}(\ODE{},...);
\end{Matlab}
will generate a variable \var{sol} which can be used with the \code{devalexp} function via
    \begin{Matlab}{OptionDOGenerator2}%
	[\var{y},\var{dy}]=\code{devalexp}(\var{sol},\var{t});
\end{Matlab}
with a vector of times \var{t}. You will get the numerical solution and its first
derivative evaluated at the times in \var{t}.

%% file: options/Refine.tex
\optshortdesc{Output refinement factor}\\
Permitted type: \refOptType{positive} \refOptType{integer} \refOptType{scalar} \code{\{1\} }

This property increases the number of output points by the specified
factor producing smoother output. Refine defaults to 1. Refine does not
apply if $\length (\var{tspan}) > 2$. To use this feature, you need a dense
output generator.

\optseealso \refOption{DOGenerator}

%% file: options/OutputFcn.tex
\optshortdesc{Output function called after each time step}\\
Permitted type and value: \refOptType{function\_handle} or \code{\{\optVal{off}\}}

If a function handle is given, it is called after each time step.
The output function has to understand the following calls:
    \begin{Matlab}{OptionOutputFcn0}%
	\code{outputFunction}([\var{t_0},\var{tfinal}],\var{y},\str{init},\var{varargin})\\%
    \code{outputFunction}(\var{t},\var{y},\str{},\var{varargin})
\end{Matlab}
The first call should initialize the output function.
The second call is executed at every time step -- either the ones selected by
the solver's error estimator, at the refined steps or at the steps
definded by \var{tspan} argument at the solver's call.
\var{t_0}, \var{tfinal} and \var{t} are scalar times. \var{t_0} and \var{tfinal}
limit the integration interval and \var{t} is the time of the current step.
\var{y} is the solution at \var{t_0} on the \str{init} call and the solution at 
\var{t} on init=\str{}. \var{varargin} is an argument which will be forwarded from
the solver call.

\optseealso \refOption{OutputSel}

%% file: options/OutputSel.tex
\optshortdesc{Indices of the solution given to the output function
}\\
Permitted type: \refOptType{vector} \refOptType{of} \refOptType{indices} \code{\{[]\} }

Only used if an output function is used. If this argument is an empty
vector, all indices of the solution are passed to the output function.
Otherwise, \var{y}(OutputSel) is passed. $OutputSel=1:\length (\var{y})$ will give the
same result as OutputSel = [].

\optseealso \refOption{OutputFcn}

%% file: options/Stats.tex
\optshortdesc{Display status messages}\\
Permitted values: \code{\optVal{silent}}, \code{\optVal{off}}, \code{\{\optVal{on}\}} or \code{\optVal{verbose}}

Set the level of status messages printed by the integrator. If \refOption{JacobianStats},
\refOption{StepStats} and \refOption{MatrixFunctionLog} are set to \str{auto}, they will be
activated only when \refOption{Stats} is set to \str{verbose} and be deactivated on all
other settings. Setting this option to \str{silent} will also disable warnings.

\optseealso \refOption{JacobianStats}, \refOption{StepStats} and \refOption{MatrixFunctionStats}

%% file: options/MatrixFunctionStats.tex
\optshortdesc{Display status messages for matrix function evaluations}\\
Permitted values: \code{\optVal{off}}, \code{\optVal{on}} or \code{\{\optVal{auto}\}}

Control, whether or not to display  status messages for the matrix
function evaluations. When set to auto, this will be activated when
\refOption{Stats} is set to \str{verbose}.

\optseealso \refOption{Stats}

%% file: options/Waitbar.tex
\optshortdesc{Show a waitbar to display progress}\\
Permitted values: \code{\{\optVal{off}\}}, \code{\optVal{on}}, \code{\optVal{text}} or \code{\optVal{both}}

Set to \str{on} to display the current progress after each time step.
Set to \str{text} if you want to display the progess on the MATLAB
prompt. Set to \str{both} if you want both kind of displays.

\optseealso \refOption{Stats}

%% file: options/ClearInternalData.tex
\optshortdesc{Clear internal integrator data}\\
Permitted values: \code{\optVal{off}} or \code{\{\optVal{on}\}}

Set to \str{off} if you want to keep the integrator's internal data.
This data is available in the global \code{eD} variable.

%% file: options/LinOp.tex
\optshortdesc{Evaluation function for the ODE's right hand side's linear part
}\\
Permitted type and values: \refOptType{function\_handle}, \code{\optVal{off}}, \code{\{\optVal{on}\}} or \refOptType{matrix}

The integrator needs to compute matrix functions with the linear part
of the right hand side of the ODE. Therefore, either the linear part has to be
evaluated, or the evaluation of the linear part times vector has to be provided --
this only works when using Krylov approximations to the matrix functions.
If this option is set to
\str{on} (default), the \ODE{} function will be called with
    \begin{Matlab}{OptionLinOp0}%
	\var{lin}=\ODE{}(\var{t},\var{y},\str{linop},{\var{varargin}});
\end{Matlab}
If it is a function handle, the call
    \begin{Matlab}{OptionLinOp1}%
	\var{lp}=\var{opts}.LinOp; \var{lin}=\var{lp}(\var{t},\var{y},{\var{varargin}});
\end{Matlab}
will be executed. Note that the linear part will only be evaluated once.

\optseealso \refOption{LinOpV} and \refOption{MatrixFunctions}

\aliases \refOption{Jacobian}

%% file: options/LinOpV.tex
\optshortdesc{Evaluation function for the ODE's right hand side's linear part
times vector}\\
Permitted type and values: \refOptType{function\_handle}, \code{\{\optVal{off}\}} or \code{\optVal{on}}

When using Krylov approximations to the matrix functions of the linear part
of the right hand side times vectors, this option can be used to provide
the result of the computation of the linear part multiplied by a vector.
Sometimes it is faster to compute this directly instead of computing the
full linear part. If this option is set to \str{on} (default), the ODE function
will be called with
    \begin{Matlab}{OptionLinOpV0}%
	\var{res}=\ODE{}(\var{t},\var{y},\str{linpop\_v},\var{v},{\var{varargin}});
\end{Matlab}
If it is a function handle, the call
    \begin{Matlab}{OptionLinOpV1}%
	\var{lpv}=\var{opts}.LinOpV; \var{res}=\var{lpv}(\var{t},\var{y},\var{v},{\var{varargin}});
\end{Matlab}
will be executed.

\optseealso \refOption{LinOp} and \refOption{MatrixFunctions}

\aliases \refOption{JacobianV}

%% file: options/LinOpStats.tex
\optshortdesc{Display status messages for the evaluation of the linear part}\\
Permitted values: \code{\optVal{off}}, \code{\optVal{on}} or \code{\{\optVal{auto}\}}

Control, whether or not to display status messages for the evaluation of the linear part.
When set to auto, this will be activated when \refOption{Stats} is set to \str{verbose}.
If enabled, this options prints
    \begin{Matlab}{OptionLinOpStats0}%
	Evaluating linear part ... done
\end{Matlab}
when the linear part is evaluated. If option \refOption{Stats} set to \str{verbose}
the Eigenvalues of the linear part will be plotted additionally.

\optseealso \refOption{Stats}

\aliases \refOption{JacobianStats}

%% file: options/Jacobian.tex
\optshortdesc{Evaluation function for the ODE's right hand side's Jacobian
}\\
Permitted type and values: \refOptType{function\_handle}, \code{\optVal{off}}, \code{\{\optVal{on}\}} or \refOptType{matrix}

The integrator needs to compute matrix functions with the Jacobian
of the right hand side of the ODE. Therefore, either the Jacobian has to be
evaluated, or the evaluation of the Jacobian times vector has to be provided --
this only works when using Krylov approximations to the matrix functions.
If this option is set to
\str{on} (default), the \ODE{} function will be called with
    \begin{Matlab}{OptionJacobian0}%
	\var{j}=\ODE{}(\var{t},\var{y},\str{jacobian},{\var{varargin}});
\end{Matlab}
If it is a function handle, the call
    \begin{Matlab}{OptionJacobian1}%
	\var{jac}=\var{opts}.Jacobian; \var{j}=\var{jac}(\var{t},\var{y},{\var{varargin}});
\end{Matlab}
will be executed.
If the Jacobian is constant, set this option to the constant matrix.

\optseealso \refOption{JacobianV} and \refOption{MatrixFunctions}

\aliases \refOption{LinOp}

%% file: options/JacobianV.tex
\optshortdesc{Evaluation function for the ODE's right hand side's Jacobian
times vector}\\
Permitted type and values: \refOptType{function\_handle}, \code{\{\optVal{off}\}} or \code{\optVal{on}}

When using Krylov approximations to the matrix functions of the Jacobian
of the right hand side times vectors, this option can be used to provide
the result of the computation of the Jacobian multiplied by a vector.
Sometimes it is faster to compute this directly instead of computing the
full Jacobian. If this option is set to \str{on} (default), the ODE function
will be called with
    \begin{Matlab}{OptionJacobianV0}%
	\var{res}=\ODE{}(\var{t},\var{y},\str{jacobian\_v},\var{v},{\var{varargin}});
\end{Matlab}
If it is a function handle, the call
    \begin{Matlab}{OptionJacobianV1}%
	\var{jacv}=\var{opts}.JacobianV; \var{res}=\var{jacv}(\var{t},\var{y},\var{v},{\var{varargin}});
\end{Matlab}
will be executed.

\optseealso \refOption{Jacobian} and \refOption{MatrixFunctions}

\aliases \refOption{LinOpV}

%% file: options/GJacobian.tex
\optshortdesc{Evaluation of the Jacobian of the nonlinear part of the ODE, for
semilinear problems}\\
Permitted type and values: \refOptType{function\_handle}, \code{\optVal{off}} or \code{\{\optVal{on}\}}

This option will only be used if \refOption{Semilin} is \str{on}.
Set this to \str{on} if the ODE file can evaluate the Jacobian of the
nonlinear part when used with flag \str{dg\_dy}. Set to a function handle if
an external function can evaluate it. Setting this to \str{off}, you will need
to use Krylov approximations to the matrix functions and set the \refOption{GJacobianV}
to something else than \str{off}.

\optseealso \refOption{Semilin}, \refOption{GJacobianV} and \refOption{MatrixFunctions}

%% file: options/GJacobianV.tex
\optshortdesc{Evaluation of the Jacobian of the nonlinear part of the ODE times
vector, for semilinear problems}\\
Permitted type and values: \refOptType{function\_handle}, \code{\{\optVal{off}\}} or \code{\optVal{on}}

This option will only be used if \refOption{Semilin} is \str{on} and MatrixFunctions
is not \str{direct}.
Set this to \str{on} if the ODE file can evaluate the product of the Jacobian
of the nonlinear part with a vector when used with flag \str{dg\_dy\_v}. Set to a
function handle if an external function can evaluate it. Setting to
\str{off} will use \refOption{GJacobian}.

\optseealso \refOption{Semilin}, \refOption{GJacobian} and \refOption{MatrixFunctions}

%% file: options/Semilin.tex
\optshortdesc{Specifies whether the ODE is semilinear}\\
Permitted values: \code{\{\optVal{off}\}} or \code{\optVal{on}}

Set to \str{on}, if the differential equation has the form
    \[y' = A y + g(t, y)\]
Then add the flags \str{gfun} and \str{dg\_dy} to your ODE file to evaluate
the non--linear part $g$ and its derivative with respect to $y$. Alternatively,
you can use the \refOption{GFun} and \refOption{GJacobian} options to set these. If you use
Krylov approximations to the matrix exponentials, you can also add the flag
\str{dg\_dy\_v} to evaluate the Jacobian of $g$ times a vector or use the
\refOption{GJacobianV} option to do so.

\optseealso \refOption{GFcn}, \refOption{GJacobian}, \refOption{GJacobianV}, \refOption{MatrixFunctions} and \refOption{JacobianV}

%% file: options/JacobianStats.tex
\optshortdesc{Display status messages for the evaluation of the Jacobian}\\
Permitted values: \code{\optVal{off}}, \code{\optVal{on}} or \code{\{\optVal{auto}\}}

Control, whether or not to display status messages for the evaluation of the Jacobian.
When set to auto, this will be activated when \refOption{Stats} is set to \str{verbose}.
If enabled, this options prints
    \begin{Matlab}{OptionJacobianStats0}%
	Evaluating Jacobian ... done
\end{Matlab}
each time the Jacobian is evaluated. If option \refOption{Stats} set to \str{verbose}
the Eigenvalues of the Jacobian will be plotted additionally.

\optseealso \refOption{Stats}

%% file: options/StepSize.tex
\optshortdesc{Step size to use}\\
Permitted type: \refOptType{non-negative} \refOptType{scalar} \code{\{0\} }

Stepsize to use. If set to $0$ (default), the integrator will choose a default
step size of $(\var{tfinal} - \var{t_0}) / 100$.

\aliases \refOption{InitialStep}

%% file: options/hConstant.tex
\optshortdesc{Use constant step size}\\
Permitted values: \code{\{\optVal{off}\}} or \code{\optVal{on}}

Uses constant step size in the integrator. In case \refOption{hConstant} is set to
\str{on}, consider setting the \refOption{InitialStep} option to specify the stepsize.

\optseealso \refOption{InitialStep}

%% file: options/InitialStep.tex
\optshortdesc{Initial step size to use}\\
Permitted type: \refOptType{non-negative} \refOptType{scalar} \code{\{0\} }

Initial step size to use. If set to $0$, the integrator will choose a default
step size of $(\var{tfinal} - \var{t_0}) / 100$. If \refOption{hConstant} is \str{on} then
the given step size here will be used for the entire integration process.

\optseealso \refOption{hConstant}, \refOption{MinStep} and \refOption{MaxStep}

\aliases \refOption{StepSize}

%% file: options/MaxStep.tex
\optshortdesc{Maximal step size to use}\\
Permitted type: \refOptType{non-negative} \refOptType{scalar} \code{\{0\} }

Set the maximal step size to use in the integration process here. Setting
\refOption{MaxStep} to $0$ (default) will result in using $(\var{tfinal} - \var{t_0}) / 10$,
so using at least ten integration steps until the integrator is finished.
This option will be ignored when using constant step size.

\optseealso \refOption{MinStep}, \refOption{InitialStep} and \refOption{hConstant}

%% file: options/MinStep.tex
\optshortdesc{Minimal step size to use}\\
Permitted type: \refOptType{non-negative} \refOptType{scalar} \code{\{0\} }

Set the minimal step size to use in the integration process here. Setting
\refOption{MinStep} to $0$ (default) will result in using $\code{eps}(\var{t})$ at integration
time \var{t}, so that $\var{t} + \var{h}$ is at least different from \var{t}.
This option will be ignored when using constant step size.

\optseealso \refOption{MaxStep}, \refOption{InitialStep} and \refOption{hConstant}

%% file: options/StepStats.tex
\optshortdesc{Display status messages for step size related events}\\
Permitted values: \code{\optVal{off}}, \code{\optVal{on}} or \code{\{\optVal{auto}\}}

Control, whether or not to display step size related status messages.
When set to auto, this will be activated when \refOption{Stats} is set to \str{verbose}.
You will be informed about step size reductions and step rejections
if this option is activated.

\optseealso \refOption{Stats}

%% file: options/exprkScheme.tex
\optshortdesc{Select Runge--Kutte scheme}\\
Permitted values: \code{\optVal{Euler}}, \code{\optVal{StrehmelWeinerA}}, \code{\optVal{StrehmelWeinerB}}, \code{\optVal{HeunA}}, \code{\optVal{HeunB}}, \code{\optVal{CoxMatthews}}, \code{\{\optVal{Krogstad}\}}, \code{\optVal{HochbruckOstermann}} or \refOptType{struct}

Select the exponential Runge--Kutta coefficients here.
The exponential Euler has no internal stage and is first--order convergent.
StrehmelWeinerA and StrehmelWeinerB are of second order, where StrehmelWeinerA
is generally preferable due to some order--loss on certain properties
of the nonlinear part of the differential equation.
HeunA and HeunB are both of order three.
CoxMatthews, Krogstad are schemes of classical order four, where HochbruckOstermann
has full order four.
You can also supply you own scheme, generate one with the
\code{rkscheme} helper function from the exprk subdirectory.

\optseealso \refOption{StrehmelWeinerParameter}

%% file: options/exprkParameters.tex
\optshortdesc{Select the free parameter(s) for some of the Runge--Kutta schemes}\\
Permitted type: \refOptType{scalar} or \refOptType{vector} \code{\{-1\} }

Select the free parameter in the StrehmelWeinerA, StrehmelWeinerB,
HeunA and HeunB schemes. In particular, the first parameter is always the
second node for the scheme. In case of HeunB, the second parameter
is the additional free parameter gamma.

%% file: options/exprbOrder.tex
\optshortdesc{Integrator order to use}\\
Permitted values: \code{\optVal{two}} (2), \code{\optVal{three}} (3) or \code{\{\optVal{four}\}} (4)

Order to use for the integrator.
When setting to two, the exponential Rosenbrock--Euler method will be used. 
The \code{\EXPRB{}32} method is the order three method with two inner
stages using the second order exponential Rosenbrock--Euler method as error estimator.
For order four, the \code{\EXPRB{}43} method will be used. This one has
three inner stages and uses a three staged order three error estimator.

%% file: options/exprbErrorEstimate.tex
\optshortdesc{Error Estimator Parameters}\\
Permitted type: \refOptType{vector} \code{\{[ 0 -2 ]\} }

The order four integrator uses an error estimator that can be tweaked with
two parameters \var{a} and \var{b}. These two parameters can be specified here with
\refOption{ErrorEstimate} = \var{[a, b]}. Choose $\var{b} \ensuremath\NEQ  1 - 6\var{a}$, otherwise the error
estimator gets close to weak order four instead of three.

%% file: options/expmssemikStep.tex
\optshortdesc{Number of old time steps to use}\\
Permitted values: \code{\optVal{one}} (1), \code{\optVal{two}} (2), \code{\optVal{three}} (3), \code{\{\optVal{four}\}} (4), \code{\optVal{five}} (5) or \code{\optVal{six}} (6)

Number of old time steps the method should use. One gives the exponential Euler
method, which is first order convergent. The number of steps is automatically
the order of the scheme as well.

%% file: options/expmssemiStartupSteps.tex
\optshortdesc{Source of the startup steps needed for the multi step method}\\
Permitted values: \code{\{\optVal{Fixpoint}\}}, \code{\optVal{Exact}} or \code{\optVal{ExpRK}}

Select the computation method for the startup steps for the multistep method.
\str{Fixpoint} uses a fixedpoint iteration to compute the steps.
\str{Exact} queries the ODE with the flag \str{exact}.
\str{ExpRK} uses an exponential Runge--Kutta scheme of the appropriate
order. Note that only methods up to order five are possible this way.

%% file: options/expmskStep.tex
\optshortdesc{Number of old time steps to use}\\
Permitted values: \code{\optVal{one}} (1), \code{\optVal{two}} (2), \code{\optVal{three}} (3), \code{\{\optVal{four}\}} (4), \code{\optVal{five}} (5) or \code{\optVal{Tokman}} (6)

Number of old time steps the method should use. One gives the exponential
Rosenbrock--Euler method, which is second order convergent. The number of
steps plus one is automatically the order of the scheme as well.
The \str{Tokman} two step scheme was proposed earlier than the more general
construction by Hochbruck and Ostermann and is third order convergent.

%% file: options/expmsStartupSteps.tex
\optshortdesc{Source of the startup steps needed for the multi step method}\\
Permitted values: \code{\{\optVal{Fixpoint}\}}, \code{\optVal{Exact}} or \code{\optVal{ExpRB}}

Select the computation method for the startup steps for the multistep method.
\str{Fixpoint} uses a fixedpoint iteration to compute the steps.
\str{Exact} queries the ODE with the flag \str{exact}.
\str{ExpRB} uses an exponential Rosenbrock--type scheme of the appropriate
order. Note that only methods up to order five are possible this way.

%% file: options/exp4DOGenerator.tex
\optshortdesc{Dense Output generator}\\
Permitted type and values: \refOptType{function\_handle}, \code{\{\optVal{exp4}\}}, \code{\optVal{hermite}} or \code{\optVal{none}}

Set the method to generate dense output. Default is \str{exp4}, \EXP4's own dense output
formula. It can be overriden to use hermite interpolation by setting this option
to \str{hermite}. Hermite interpolation is not applicable for stiff problems.
Use it for testing only. Additionally the dense output generator can be switched off.

%% file: writingMatrixFunctions.tex
\ifdefined\mathe
	\def\globalJac{\code{\var{eD}.\var{jac}.\var{J}}\xspace}
	\def\globalGJac{\code{\var{eD}.\var{jac}.\var{Jg}}\xspace}
	\def\theGlobal{}
\else
	\def\globalJac{\code{\var{jac}}\xspace}
	\def\globalGJac{\code{\var{gjac}}\xspace}
	\def\theGlobal{the global\xspace}
\fi

\Abschnitt{Matrix Functions}{matFun}
	This section deals with custom evaluation functions for the product of matrix functions with vectors.

	\INTCMD has two built-in methods to compute these evaluations: directly by diagonalisation and using a Krylov subspace method. In some situations, it can be useful to provide a custom function that does this job. If the Jacobian has a special structure that makes it possible to compute matrix exponentials differently even in large dimensions, this would be such a case. For this reason, it is possible to hook an arbitrary function into the integration process. See the \refOption{MatrixFunctions} option on how to accomplish that.
	
	The two internal functions are implemented the same way as an external method would have to work. They can be used for reference and can be found in the \EXPODE distribution at \file{matFun/matFunDirect.m} and \file{matFun/matFunKrylov.m}. The former one is much easier to understand due to its size and complexity.
	
	It will be neccessary to use the global \var{eD}%
	\ifdefined\mathe%
		variable which contains
	\else%
		, \var{jac} and \var{gjac} variables which contain
	\fi
	information shared between the various \EXPODE functions. See subsection \refUnterabschnitt{eD} for details on
	\ifdefined\mathe
		this variable.
	\else
		these variables.
	\fi
	The parts important for this context will be explained when they appear.
	
	All logging output of the function should be printed with the \code{\var{eD}.\var{log}.\var{matFunLog}} function. This way it can easily be redirected or disabled by the integrator.
	
	The custom function will be called \var{matFun} in the remainder of this section. It will be invoked in different ways by the \EXPODE integrators. The function's signature needs to have the following form:
	
% 	function [h, varargout] = matFunDirect(job, t, y, h, flag, v, reusable, reuse, facs)
	\begin{Matlabfun}{matfunhead}%
		\FUNCTION [ \v{h}, \v{varargout} ] = \ifx\mathe\undefined...\\%
		~~~~\fi\v{matFun}(\v{job}, \v{t}, \v{y}, \v{h}, \v{flag}, \v{v}, \v{reusable}, \v{reuse}\ifx\mathe\undefined, \v{facs}\fi) \\%
			\GLOBAL eD; \\%
			\v{o} = \v{eD}.\v{int}.\v{o}; \\%
			\v{funs} = \v{eD}.\v{functions};
	\end{Matlabfun}
	
	The \code{\var{flag}} variable has two different roles. It can either request a specific action to be performed or it can inform the method about the meaning of the currently evaluated matrix function. The request-type flags will be discussed first.
	\begin{Matlab}{matfuninit}%
		\v{matFun}([], [], [], [], \s{init});
	\end{Matlab}
	should initialize the function globally. It has to state whether or not it requires to explicitly evaluate the Jacobian of the right hand side and -- if used -- the one of the non-linear part $g$ for semilinear problems (see option \refOption{Semilin}). This is done by setting
	\begin{Matlabfun}{matfunJacExplicit}%
		\v{eD}.\v{matFun}.\v{needJacExplicit} = true; \\%
		\v{eD}.\v{matFun}.\v{needGJacExplicit} = true;
	\end{Matlabfun}
	or to \code{false} respectively. The direct solver for instance sets both values to \code{true}, since it diagonalizes the Jacobian, the Krylov method sets it so \code{false}, since it only needs to evaluate the product of the Jacobian with a vector. If the first variable is \code{true}, then \theGlobal \globalJac variable contains the evaluation of the Jacobian at the current step. \globalGJac contains the evaluation of $g$'s Jacobian.
	
	Additionally, all statistical data fields have to be initialized. They are later used with the \str{statistics} flag. Example: the direct evaluator sets
	\begin{Matlabfun}{matfunStatsInit}%
		\v{eD}.\v{stats}.\v{matFun}.\v{NofDiag} = \l{0}; \\%
		\v{eD}.\v{stats}.\v{matFun}.\v{NofMFEv} = \l{0};%
	\end{Matlabfun}
	These two variables count the number of matrix function evaluations and the number of diagonalizations of the Jacobian.
	
	We continue with the registerjobs phase. It will be called, after the solver step decides which matrix functions it needs to evaluate. This can be in two situations: Either after the matrix function evaluator was initialized with the \str{init} flag and now has enough information about the differential equation and the options, or when the integration scheme requires a change of the matrix function coefficients or even exchanges some matrix functions. The latter can be the case for instance in the multistep integrators, after their startup phase is finished and the normal integration routine is started. The registerjobs phase is triggered by
	\begin{Matlab}{matfunregisterjobs}%
		\v{matFun}(\v{jobs}, [], [], [], \s{registerjobs});
	\end{Matlab}
	
	The \var{jobs} variable is a struct of vectors or matrices. The fieldnames of this struct represent the meaning of the vector to multiply the matrix function with. Typically \str{F} is used for a product with the right hand side, and \str{v} for the non-autonomous correction. Let \var{job} be a field of \var{jobs}. Each row of the matrix \var{job} contains the coefficients for a linear combination of the matrix functions in \code{\var{eD}.\var{int}.\var{jobFunctions}}. \code{\var{eD}.\var{int}.\var{jobFunctions}} contains the scalar, vectorized versions of the matrix functions. For \EXPRB it is constructed the following way:
	\begin{Matlabfun}{matfunJobsfunctions}%
		\v{eD}.\v{int}.\v{jobFunctions} = \{\ind \\%
			@\v{phi1} \\%
			@\v{phi2} \\%
			@\v{phi3} \\%
			@\v{phi4} \\%
			@(\v{A}, \v{h}) \v{phi1}(\v{A}, \v{h}/\l{2}) \\%
			@(\v{A}, \v{h}) \v{phi2}(\v{A}, \v{h}/\l{2})\unind \\%
		\};
	\end{Matlabfun}
	The above functions are based on the \var{phim} function from the \code{EXP4} package \cite{exp4}. If \var{job} is an element of \var{jobs} for line number \var{k} in a \var{job} matrix, \var{matFun} should return the following in the \var{k}th column of the result variable \var{res} (written in a suggestive notation):
	\begin{Matlabfun}{matfunJobOutput}%
		\v{res}\{(:,\v{k})\} = $\sum\limits_{\var{m} = 1}^{\command{length}(\var{job}\{\var{k}\})}$ \v{job}\{\v k\}(\v m) * (\v{eD}.\v{int}.\v{jobFunctions}\{\v{m}\}($\code{J}_{n}$, h) * v);%
	\end{Matlabfun}
	
	The next flag is the initstep flag:
	\begin{Matlab}{matfuninitstep}%
		\v{matFun}([], \v{t}, \v{y}, \v{h}, \s{initstep});%
	\end{Matlab}
	When this flag is given, the function should prepare itself for several matrix function evaluations with the same \var{t}, \var{y} and \var{h}.
% 	This way, \var{matFun} knows the functions it has to evaluate. 
	%The functions \var{phim}, \var{psim}, \var{phi3m} and \var{phi4m} are part of the \code{EXP4} package \cite{exp4} and are included in \EXPODE.
	
	At the end of the integration process, \var{matFun} will be called with flag \str{cleanup}:
	\begin{Matlab}{matfunCleanup}%
		\v{matFun}([], [], [], [], \s{cleanup});%
	\end{Matlab}
	Then it should clean up all persistent variables it uses. The \var{matFun} field of the global \var{eD} variable will automatically be cleared by the integrator.
	
	The call
	\begin{Matlab}{matfunDesc}%
		\v{desc} = \v{matFun}([], [], [], [], \s{description});%
	\end{Matlab}
	should return a basic information string, naming the method used to calculate the results. The result will be printed at the beginning of the integration phase by
	\begin{Matlab}{matfunDescPrint}%
		\command{sprintf}(\str{Matrix functions evaluated \%s.$\backslash{}$n}, ... \\%
		\ind\ind\indent matFun([], [], [], [], \str{description}));
	\end{Matlab}
	
	The last flag, which is not used for a computational purpose, is
% 	The last that does not request an actual computation is
	\begin{Matlab}{matfunStats}%
		\v{matFun}([], [], [], [], \s{statistics});%
	\end{Matlab}
	Then \var{matFun} is expected to print some statistics. See the \str{init} flag above for an example on what the direct method logs. If there are no interesting statistical data collected, this flag should be ignored.
	
	All other flags that are provided should be used by \var{matFun} to recognize the use of the executed job. The call looks like this:
	\begin{Matlab}{matfuneval}%
		\v{matFun}([], \v{t}, \v{y}, \v{h}, \v{flag}, \v{v}, \v{reusable}, \v{reuse}\ifx\mathe\undefined, \v{facs}\fi));
	\end{Matlab}
	\var{flag} needs to be of the fields of the \var{jobs} structure provided in the registerjobs phase. In that case, \var{matFun} needs to actually compute the result of the matrix functions specified by \var{job} evaluated at \code{h * $\code{J}_{n}$} or \code{h * A} multiplied by the vector \var{v}. The required linear combinations of the $\varphi$-functions for this computation have already been supplied at the registerjobs phase above.
	
	% TODO
	
	There are two additional parameters: \var{reusable} and \var{reuse}. The first one indicates whether the matrix function evaluated at the {\itshape same} arguments will have to be used again if the current step was rejected by the error estimator. The second one will be set \code{true} if the previous step was actually rejected and the values calculated there can be reused. In case of an exact evaluation, the diagonalization computed the last time can be reused. When approximating the results (e.g. in the Krylov version) it may be required to calculate up to a higher precision. In case of Kylov approximizations, the old Krylov subspaces can be extended. To save the reusable data, use the structure object
	\begin{Matlabfun}{matfunReusable}%
		\v{eD}.\v{matFun}.\v{save}.(\v{flag})
	\end{Matlabfun}
	This is one of the reasons to provide those \var{flag}s. The other reason is to save statistical data separately for each matrix function call in the integration scheme. Use
	\begin{Matlabfun}{matfunStatisticsPerMatFun}%
		\v{eD}.\v{stats}.\v{matFun}.(\v{flag})
	\end{Matlabfun}
	as storage here.
	
	For the addition of \EXP4 to the \EXPODE package, another feature was needed, namely the evaluation of $\varphi_1(j h J_n) v$ and $\varphi_2(j h J_n) v$, $j = 1, ..., \code{\var{facs}}$. This needs to be done at once when the \code{\var{facs}} argument is given with $\code{\var{facs}} > 1$. The functions \code{\var{phi1}} and \code{\var{phi2}} were extended to be called $\code{\var{res} = \var{phi1}(\var v, \var h, \var{facs})}$ and return an $n \times 3$ Matrix containing the required evaluations at once.
	
	The implementation of custom matrix functions is quite a complex task. The steps needed were described as simple as possible even though they are still quite hard to understand only by reading this section. Therefore, it is highly recommended to look at the two existing implementations as a reference.

%% file: writingexpode.tex
\Abschnitt{Writing an \code{EXPODE} Integrator}{writingexpode}

This section contains information for integrator authors. It overviews the \EXPODE helper routines and how and where to use them. It is split up into several subsections.

The first subsection describes the fundamentals of the \EXPODE package. It names the code files needed, explains the package structure and sets coding style standards to make the \EXPODE code as a whole consistent and well readable.

The second subsection deals with the integrators themselves. It states the calling conventions which will make the usage consistent. It will introduce the helper functions that will ease the author's life. These will let him concentrate on the integrator's details instead of having to deal with memory management or syntactic correctness of options. Semantic correctness still has to be considered by the author.

The third subsection deals with integrator options. Here we will explain how to create the option description structure.

The last subsection explains the global \var{eD} variable. This variable contains all global information shared by \EXPODE functions.

\Unterabschnitt{\code{EXPODE} Basics}{expkitbasics}

\ifx\mathe\undefined
	There are two different packages for \EXPODE. The first one is the normal \EXPODE integrator package, containing all fully working integrators and all their required helper functions. This package is targeted at normal users. The second one, \EXPODEDEV is an extension of the first one. It additionally contains some tools specifically for developers of new integrators.
\fi

The \EXPODE package consists of the following parts:
\begin{itemize}
	\item the integrators (\EXPODECMD, \EXPRB, etc),
	\item the integrator info functions (\code{exprbinfo}, etc),
	\item the options helper functions (\EXPSET, \EXPRBSET, etc),
	\item the \code{initPaths} routine to setup the appropriate paths,
	\item the \DEVALEXP function to evaluate the numerical solution at arbitrary times,
	\item common helper functions for the integrators in the \path{expode} directory,
	\item general helper functions in the \path{helpers} directory,
	\item matrix-function evaluation helpers in the \path{matFun} directory,
	\item some files from the \code{exp4} package \cite{exp4} in the \path{3rdparty/exp4} directory,
	\item a set of helper functions for each integrator in folders with the names of the integrators and
	\item some examples in the \code{examples} directory.
\end{itemize}
Each integrator has to provide at least the integrator M-file, the options helper function, the information command, a function that creates the option description structure and a setup command that contains information for the basic integration routine \EXPODECMD and evaluates integrator specific options. The latter two have to be placed in the helper directory specific for the integrator. As an example see the files \file{exprb.m}, \file{exprbset.m}, \file{exprbinfo.m}, \file{exprb/exprbOpts.m} and \file{exprb/exprbSetup.m}. The integrator specific helper directory can also contain additional supportive functions which will be accessible from the integrator.

\ifx\mathe\undefined
	The \EXPODEDEV package has the following additional components:
	\begin{itemize}
		\item a function \code{generateOptionDoku} to convert the in-\MATLAB help for the options to \LaTeX ~ code for inclusion in a documentation like this one,
		\item the \BASH script \path{createRelease.sh}, which creates release packages,
		\item a directory \path{unitTests} containing tools for automated tests -- useful for finding regressions,
% 		\item another directory \path{test} containing some non-automated tests for some of the components,
		\item a last directory \path{prototype} containing stubs for the minimally required functions for a new integrator and a script \path{generateIntegrator.sh} to create a new \EXPODE integrator from them.
	\end{itemize}
\fi

For readability reasons, the code files should obey the following conventions:
\begin{itemize}
	\item All M-files should not contain any warnings (and of course no errors) displayed by \MATLAB's editor. Usually warnings either mark bad coding style or optimization possibilities. Rarely it is required to disable single warnings like unused function parameters.
	\item All names of functions, that are direcly visible to the user should be completely lower case. These are mainly the functions in the top level directory except \code{initPaths} which is used internally to setup paths. Example \code{exprbinfo}.
	\item All names of variables and internal functions should begin with a lower case letter. Each new word in the names should start with an upper case letter. Example: \code{checkValidOptions}, \var{quitIntegrator}
	\item Each M-file should contain a help text in its header. Longer functions that do non-trivial work should have comments in the code as well.
	\item Indentation is done via four spaces, no tabs.
	\item All commas should be followed by a space. Insert spaces after an open squared bracket or brace and before a closing one. Empty brackets can be written without spaces. Use no spaces around parentheses. Equal signs should be surrounded by spaces. Example:
	\begin{Matlab}{spaces}%
		\v{result} = \command{str2func}([ \s{exprb}, \s{opts}, [] ]);
	\end{Matlab}
	\item All lines except the ones that are part of control strucures like \command{if}s or \command{while}s should be terminated by a semicolon.
	\item Lines should be maximally around eighty characters log. \quot{Around} means that a line should be split up after the first word that exeeds the eightieth column. \MATLAB can highlight this column. This option can be activated in the preferences dialog under \quot{Editor/Debugger} $\rightarrow$ \quot{Display}. The \MATLAB syntax supports split lines via three dots at their ends. Strings can be split up, too:
	\begin{Matlab}{brokenString}%
		\v{text} = [ \s{This is a very very very } ... \\%
		~~~~~~~~\s{very long text} ];
	\end{Matlab}
	Note the space after the last \quot{very} in the first line, because the variable \var{text} will not contain a newline chacter where the text is split. The remainder of the split lines has to be indented eight more spaces than its parent.
	\item One-line-\command{if} constructions should be avoided as well as nested expressions. Generally it is dicouraged to put more than one command in a single line. Do not use constructions like:
	\begin{Matlab}{onelineif}%
		\IF condition, expression,~\indent \END
	\end{Matlab}
	or
	\begin{Matlab}{nested}%
		\v{var1} = func1(func2(\v{var2}, func3(\v{var3})), func4(\v{var4}));
	\end{Matlab}
	The former should be written in three lines, the latter at least expanded into
	\begin{Matlab}{splitnested}%
		\v{temp1} = func2(\v{var2}, func3(\v{var3})); \\%
		\v{temp2} = func4(\v{var4}); \\%
		\v{var1} = func1(\v{temp1}, \v{temp2});
	\end{Matlab}
	\item Function code should always be indented one level.
	\item All functions should be terminated with an \command{end}. Nested functions should be indented. Additional functions in one M-file should \emph{not} be intended. The latter ones should only be used, if it is clear, that no other function could use the additional function in the future. A file \file{majorFunction.m} should look like:
	\begin{Matlabfun}{functionIndentation}%
		\FUNCTION majorFunction \\%
			\c{majorFunctionCode} \\%
			\\%
			\FUNCTION nestedFunction \\%
				\c{nestedFunctionCode} \\%
			\END \\%
		\END \\%
		\\%
		\FUNCTION additionalFunction \\%
			\c{additionalFunctionCode} \\%
		\END
	\end{Matlabfun}
% 	\begin{Matlabfun}{functionIndentation}%
% 		\FUNCTION majorFunction \\
% 			majorFunctionCode \\
% 		\\%
% 		\unind\FUNCTION nestedFunction \\
% 			nestedFunctionCode
% 	\end{Matlabfun}
	\item There is {\itshape only one} global variable to store common data, a \command{struct} called \var{eD}, and two additional ones to store the Jacobians of the right hand side of the differential, \var{jac} and a non-linear part in case of seminilinear equations, \var{gjac}. \EXPODE functions should never define additional ones. Generally, it is discouraged to use the \var{eD} variable, since access to it is slow due to some \MATLAB restrictions. Use it only if really necessary or the information stored are of common interest doublessly. The two exceptions of this \emph{only one} rule -- even though unclean by design -- are acceptable since \MATLAB has a bottleneck accessing big datablocks in global variable structures. % As for each strict rule, there is one exception. There are two other global
% 	\item There is {\itshape only one} global variable: \var{eD}. \EXPODE functions should never define additional ones. Generally, it is discouraged to use this variable. Use it only if really necessary or the information stored are of common interest doublessly.
\end{itemize}

\ifdefined\mathe
	\input{writingexpodeints-maths}
\else
	\input{writingexpodeints-info}
\fi

\Unterabschnitt{Options}{expodeoptions}

This section deals with integrator options. \EXPODE has advanced option handling subroutines which allow automatic checking of options set by the user. First, each integrator needs an options helper function as a user interface. This function should be called the same as the integrator with a \quot{set} suffix. The \EXPRB integrator is contained in the \file{exprb.m} M-file and its options helper is named \quot{\path{exprbset.m}} for instance. The options helper should simply be a wrapper around the \EXPSET function which only adds the integrator name as the first argument:
\begin{Matlabfun}{exprbset}%
	\FUNCTION \v{opt} = expnewset(\v{varargin}) \\%
		\v{intname} = \s{expnew}; \\%
		initPaths(\v{intname}); \\%
		\\%
		\IF \ISEMPTY{}(\v{varargin}) \\%
			\HELP expnewset; \\%
			\RETURN; \\%
		\END \\%
		\\%
		opt = expset(intname, varargin\{:\}); \\%
	\END
\end{Matlabfun}
The \code{initPaths} routine is needed to initialize the \EXPODE internal paths as seen in the previous section.
\ifx\mathe\undefined
	This method can be automatically generated by the \code{generateIntegrator.sh} script, see the previous section.
\fi

The \code{expset} method needs a so-called \emph{option description object}. This has to be provided by a function called \code{expnewOpts} and has to be located in the \path{expnew} directory. This object contains all information about the options for \EXPNEW. It is a \MATLAB struct object and has to have the following strutcure:

\begin{Matlabfun}{optDescObject}%
	>{>} \EXPNEW{}Opts \\%
	\var{ans} = \\
	~~~~~~~~integrator: 'expnew' \\
	~~~~integratorname: 'exponential new-type integrator' \\
	~~~~integratordesc: [1x223 char] \\
	~~~~~~~~~~~~~usage: [1x522 char] \\
	~~~~~~~~~~~~~~opts: [1x1 struct] \\
\end{Matlabfun}
The first four fields are all strings, the first one containing the integrator name, the second a longer name of the integration method, the third a short description, and the last one the help text for the integrator, displayed when calling the integrator without arguments.
\ifdefined\mathe
	This is the reason the \code{initIntegrator} method can set the \var{quitIntegrator} variable to \code{true}. It checks if no arguments were provided and prints this help message.
\else
	This will be done automatically when wrapping the \EXPODECMD base integrator, see the \hrefUnterabschnitt{expodeintegrators}{previous section}.
\fi

The last entry, \var{opts}, is again a structure object and it contains the actual options. Each element of this structure represents an option where the name of the element ist the name of the option. The option itself is a struct once again and can be created via
\begin{Matlabfun}{optNewOption}%
	\var{opts}.NewOption = createOptDesc( ... \ind \\%
		\str{short description}, ... \\%
		\str{type}, ... \\%
		\str{default}, ... \\%
		\{ \str{value1}, \str{value2} \}, \{ ... \\%
		\str{Very long and detailed description for the very nice} \\%
		\str{and new option NewOption.}\\%
		\}, \{ \\%
		\str{OptionA}, \str{OptionB} \\%
	\unind\});
\end{Matlabfun}
The first argument is a short description that will be displayed when using \code{\EXPNEW{}info}, see below. The next one contains the option's type. Available option types are listed in section \refAbschnitt{expset}. Next, we have the default value for the option. This can be a numeric value or string, depending on the option type. In case of \refOptType{boolean} or \refOptType{list} it can either be the numeric or string value of the option. Argument four of the \code{createOptDesc} function contains the allowed values for \refOptType{list}. For all other types, provide the empty cell object \code{\{\}}. The next argument contains the long and detailed description for the option. This can either be a string or a cell of strings as used in the example. The last argument contains a list of referenced options. This list will be printed in the form \code{See also: OptionA, OptionB} when calling \code{expnewinfo NewOption}.

\ifx\mathe\undefined
	For advances usage two additional arguments can be given, \var{renameTo} and \var{aliases}, in that order. The first one is used at the moment, where options are evaluated to their numeric values. The option will then be renamed to \var{renameTo}. This is useful if the same option has different names in different scenarios. In case of variable step size integrators, there is an option \refOption{InitialStep} which indicates the step size to start with, where in the case of constant step size options there is the \refOption{StepSize} option, which will be renamed to \code{InitialStep} in the integration process to avoid code duplication.
\fi

For reference and a lot of examples look at the file 
\ifdefined\mathe
	\code{exprb/exprbOpts.m}
\else
	\code{expode/baseOpts.m} or any other \code{*/*Opts.m} file
\fi
in the \EXPODE package.
	
\ifx\mathe\undefined
	A basic stub for the opts command will be generated by the \code{generateIntegrator.sh} script. There are a set of basic options the integrator has to provide, see section \refAbschnitt{expodeopts} for details. They can be generated with
	\begin{Matlabfun}{optNewBaseOpts}%
		\var{optionDesc} = baseOpts('expnew', 'exponential new-type integrator');
	\end{Matlabfun}
	This line of \code{expnewOpts.m} is part of the stub and doesn't need to be adjusted. These basic options have to be extended either for use in semilinear integrators or linearized integrators and for constant or variable step size integrators. The two code lines responsible for this are
	\begin{Matlabfun}{optNewStandardOpts}%
		\var{optionDesc} = linearizationOpts(\var{optionDesc}); \\%
		\var{optionDesc} = varStepOpts(\var{optionDesc});
	\end{Matlabfun}
	Use \code{semilinOpts} instead of \code{linearizationOpts} for semilinear integrators, \code{constStepOpts} instead of \code{varStepOpts} for constant step size integration. Using \code{constStepOpts} automatically switches off the step size estimator.
\fi

% \ifx\mathe\undefined
% 	A basic stub for the opts command will be generated by the \code{generateIntegrator.sh} script.
% \fi

Another function that needs to be provided for each integrator is the \code{info} command. This is a wrapper around \code{optInfo} similar to the options setter. It will print information on the integrator. The information will be extracted from the option description object described above. Here is its basic structure:
\begin{Matlabfun}{expnewinfo}%
	\FUNCTION expnewinfo(\var{optName}) \\%
		\IF \var{nargin} < 1 \\%
			\var{optname} = \str{{}}; \\%
		\END \\%
		optInfo(\str{expnew}, optName); \\%
	\END
\end{Matlabfun}
\ifx\mathe\undefined
	Like the set, opts and the integrator command this file will be generated by the script \code{generateIntegrator.sh}.
\fi

\Unterabschnitt{The Global \code{eD} Variable}{eD}

This section deals with the main global variable \var{eD}. Except for \var{jac} and \var{jacv} this is the only global variable used throughout \EXPODE (see the coding conventions above). Since it has such basic relevance it will be discussed in detail here. It is a structure object containing several different types of information. After an
\ifdefined\mathe
	\EXPRB call
\else
	integrator call, here \EXPRB as an example,
\fi
it contains the following fields if the \refOption{ClearInternalData} is set to \str{off}:
% There is only one global variable used throughout \EXPODE (see the coding conventions above), namely \var{eD}. It is a structure object containing several different types of information. After an \EXPRB call it contains the following fields if the \refOption{ClearInternalData} is set to \str{off}:

\begin{Matlabfun}{eDContent}%
	>{>} \v{eD} \\%
	\v{eD} = \\%
	~\v{integrator}: \str{exprb} \\%
	~~~~~~~~\v{ODE}: [1x1 struct] \\%
	~~~~~~~~\v{int}: [1x1 struct] \\%
	~~~~~~~~\v{log}: [1x1 struct] \\%
	~~~~~~\v{evals}: [1x1 struct] \\%
	~~\v{functions}: [1x1 struct] \\%
	~~~~~~\v{stats}: [1x1 struct] \\%
	~~~~~\v{matFun}: [1x1 struct] \\%
	~~~~~~\v{tdata}: [451x1 double] \\%
	~~~~~~\v{ydata}: [451x19 double]
\end{Matlabfun}

The first field is \var{integrator} and it holds the name of the integrator. The second, \var{ODE}, contains information about the differential equation and the solver request:

\begin{Matlabfun}{eD.ODEContent}%
	>{>} \v{eD.ODE} \\%
	\v{ans} = \\
	~~~~~~~~\v{t0}: \literal 0 \\
	~~~~~~~~~\v{T}: \literal{45} \\
	~~\v{duration}: \literal{45} \\
	~~~~~~\v{tDir}: \literal{1} \\
	~~~~~~~~\v{y0}: [19x1 double] \\
	~~~~~~\v{ylen}: \literal{19} \\
	~\v{transpose}: \literal{1} \\
	~~~\v{options}: [1x1 struct] \\
\end{Matlabfun}

It contains the endpoints of the integration interval, \var{t0} and \var{T}, their difference \var{duration}, the integration direction (\var{tDir}) -- with values \literal{1} for increasing or \literal{-1} for decreasing time steps, the initial value \var{y0} and its length \var{ylen} as well as the \var{options} provided at the integrator's call. \var{transpose} is a boolean that indicates, whether the solution vectors are rows or columns. In the latter case, they have to be transposed before being stored in the result matrix.

The \var{int} field of \var{eD} contains information for the integration process. Is has the following form:
\begin{Matlabfun}{eD.intContent}%
	>{>} \v{eD.int} \\%
	~~~~~~~~~~~~~~~~~~~~\v{o}: [1x1 struct] \\%
	~~~~~~~~~~~~~~~~~~~\v{nv}: [1x1 struct] \\%
	~~~~~~~~~~~~~~~~~~~\v{od}: [1x1 struct] \\%
	~~~~~~~~~\v{jobFunctions}: {6x1 cell} \\%
	~~~~~~~~~~~~~~~~\v{order}: 4 \\%
	~~~~~~~~~~~\v{errorOrder}: 4 \\%
	~~~~~~~~~~~~\v{multiStep}: 1 \\%
	~~~~~~~~~~~~~~\v{semilin}: 0 \\%
	~\v{denseOutputGenerator}: @hermite \\%
	~~~~~\v{wantSolverOutput}: 0 \\%
	~~~~~~\v{wantDenseOutput}: 1 \\%
	~\v{solverOutputStepData}: \{[19x1 double]  [19x1 double]\}
\end{Matlabfun}

The fields \var{o} contains a version of the options structure but with all non-set options evaluated to their default values and all \refOptType{list}-type options evaluated to the index of the selected value in the list. \refOptType{boolean}-type options are also converted to \literal{1} for value \code{true} and \literal{0} for \code{false}. The last field, \var{nv}, contains the numeric values of the list-type options, see the previous section. The field \var{od} contains the option description strcture describing the integrator's options that were used to evaluate \var{eD.ODE.options} into \var{eD.int.o}. \var{jobFunctions} contains function handles to scalar, vectorized versions of the matrix functions involved in the integrator's scheme. These will be used by matrix function evaluators, see sections \refAbschnitt{matFun} and and \refUnterabschnitt{expodeintegrators}. The field \var{order} represents the integrator order, whereas \var{errorOrder} ist the order of the error estimator plus one -- i.e. the parameter for the \code{errstep} function used to compute the next timestep size. The entry \var{multiStep} contains the number of old time step data necessary to compute the next time step, \var{semilin} indicates, whether the solver is for semilinear problems, i.e. whether in needs to evaluate the linear part once at the beginning or the Jacobian at each time step. The \var{denseOutputGenerator} field is a handle to the dense output generation function. The two flags \var{wantSolverOutput} and \var{wantDenseOutput} indicate, wheter the user requested one of the two outout types. If the first is \literal{true}, then \var{solverOutputStepData} contains the data -- additionally to the time step data -- needed to evaluate the numerical solution by the \DEVALEXP function. If the second is \literal{true} it contains the data necessary for the dense output evaluation between the last two time steps.

% The \var{int} field of \var{eD} contains three entries: \var{jobFunctions} contains function handles to scalar, vectorized versions of the matrix functions involved in the integrator's scheme. These will be used by matrix function evaluators, see sections \refAbschnitt{matFun} and and \refUnterabschnitt{expodeintegrators}. Additionally, there exist the fields \var{o} containing a version of the options structure but with all non-set options evaluated to their default values and all \refOptType{list}-type options evaluated to the index of the selected value in the list. \refOptType{boolean}-type options are also converted to \literal{1} for value \code{true} and \literal{0} for \code{false}. The last field, \var{nv}, contains the numeric values of the list-type options, see the previous section.  %Example: There is the \refOption{MatrixFunctions} option, allowing the values \str{direct} and \str{arnoldi} (in that order). Then \var{o}.\var{MatrixFunctions} is \literal{0} for a direct solution and \literal{1} for the use of the Arnoldi method. Then we have fields \code{\var{nv}.\var{MatrixFunctions}.\var{direct} = 0} and \code{\var{nv}.\var{MatrixFunctions}.\var{arnoldi} = 1}. So somewhere in the integrator code one could find something like
% \begin{Matlabfun}{eD.int.nv}%
% 	\SWITCH \v{eD}.\v{int}.\v{o}.\v{MatrixFunctions} \\
% 		\CASE \v{eD}.\v{int}.\v{nv}\var{MatrixFunctions}.\var{direct} \\
% 			\c{do something} \\
% 		\CASE \v{eD}.\v{int}.\v{nv}\var{MatrixFunctions}.\var{arnoldi} \\
% 			\c{do something else}
% \end{Matlabfun}

Next in the \var{eD} structure we have the \var{log} field containing logging functions for several purposes. They all have the same syntax as the \command{fprintf} command, but may be directed to a file or simply ignored depending on what the user wants to see.
\begin{Matlabfun}{eD.logContent}%
	>{>}  eD.log\\%
	\v{ans} = \\%
	~~~~~~\v{verbose}: @fprintf \\%
	~~~~~~~\v{status}: @(varargin)wrapped(fun,prefix,varargin\{:\}) \\%
	~~~\v{statistics}: @fprintf \\%
	~~~~~~~\v{jacLog}: @nullfunction \\%
	~~~~~~\v{stepLog}: @fprintf \\%
	~~~~\v{matFunLog}: @nullfunction \\%
	~~~~~~\v{warning}: @warning \\%
	~~~~~~~~\v{error}: @error
\end{Matlabfun}

In \var{eD.functions} we can find evaluation functions for the ODE and its Jacobian, furthormore we find a handle to the matrix function evauator and last we can retrieve the \var{varargin} parameters provided by the user at the integrator's call. In the example the field reads
\begin{Matlabfun}{eD.functionsContent}%
	>{>}  eD.functions\\%
	\v{ans} = \\%
	~~~~~~\v{ode}: @(t,y,varargin)ode(t,y,flag,varargin\{:\}) \\%
	~~~~~~\v{jac}: @(t,y,varargin)ode(t,y,flag,varargin\{:\}) \\%
	~~\v{matFunV}: @matFunKrylov \\%
	~\v{varargin}: \{\}
\end{Matlabfun}

The \var{stats} field in \var{eD} contains statistical data collected over the whole integration process. It has a subfield \var{matFun} which contains data gathered by the matrix function evaluators.

The fields \var{tdata} and \var{ydata} contain the data of the numerical solution. This is what the integrator returns to the user.

% The last field, \var{evals} contains evaluations of the right hand side of the differential equation, its Jacobian and the Jacobian of the non-linear part -- if \refOption{Semilin} is \str{on}.

Most of this data will be cleaned before the integration ends. If you want to preserve the data for testing, set the option \refOption{ClearInternalData} to \str{off}.

\ifx\mathe\undefindex
	\Unterabschnitt{Deploying \EXPODE}{deployment}
	The previous subsection concludes the basic information needed to write new integrators. Additional help can be found using \MATLAB's \command{help} command on the \EXPODE functions. All should be documented well both helptext and code wise. In this subsection you will find some advices for deployment of your newly created integrator.
	
	The \EXPODEDEV package contains a unit testing framework. Before packaging the new release it is highly recommended that you run at least the existing tests to make sure that you didn't introduce any regression to the previous release. This is done by
	\begin{Matlab}{unitTests}%
		cd \literal{unitTests}; \\%
		runAllTests;
	\end{Matlab}
	If all goes well you should see a message like
	\begin{Matlabfun}{unitTestResults}%
		Expode unit tests framework finished. \\%
		Totally ran 593 tests, out of where 593 finished successful and 0 failed. \\%
		Overall testing time: 00:01:37.
	\end{Matlabfun}
	Additionally you should write new unit tests for you new components und run them. Tests have to be located in the \path{unitTests/tests} directory. The testing framework clears all global variables and restores \MATLAB's default path before calling the actual test. Only one global variable for the testing data will be created, the test itself does not have to use it, though, there are helper functions for it. The test's working directory is \path{unitTests/tests}. You will need to call an intialization routine to tell the framework what you want to test. The next thing you need to do is add the path to your tested functions.
	\begin{Matlabfun}{ownUnittest}%
		startTest(\str{optionDesc}); \\%
		addpath ../../expode;
	\end{Matlabfun}
	Now you can write your testing code. You may need set some values in the global \var{eD} variable. After you are done, use
	\begin{Matlabfun}{ownUnittestEnd}%
		finishTest(\v{success});
	\end{Matlabfun}
	to report the result, where \var{success} is a boolean containing your result. Running this specific test is done via
	\begin{Matlab}{unitTestsRunSingle}%
		runTest \literal{testname};
	\end{Matlab}
	run from the \path{unitTests} directory again, where \code{testname.m} is the name of the m-file. Note that you can include more that one test into one file, you just have to separate them by \code{startTest} and \code{finishTest}. The \code{runTest} routine can also run multiple tests in a row, just give more testnames separated by a space.
	
	Now that you finished testing you new integrator, it is time to release it. In \EXPODEDEV you will find a script, \code{createRelease.sh} to do that. You have to add the integrator to the \code{ints} variable. This will automatically include the integrator m-file, the options helper, the info command and the integrator's directory. If you have additional files, you have to add it to the apropiate variable in the script. This should usually not be the case, since all methods that are only used by your new integrator should be in its own directory. Now you can call
	\begin{Matlabfun}{bash:releaseInt}%
		/path/to/expode \$ ./createRelease.sh 1.1
	\end{Matlabfun}
	to create a package for \EXPODE, version 1.1. To generate a \EXPODEDEV package, call
	\begin{Matlabfun}{bash:releaseIntDev}%
		/path/to/expode \$ ./createRelease.sh 1.1-dev -d
	\end{Matlabfun}
	which will result in version 1.1-dev with the additional functions described at the beginning of this section.
\fi

%% file: writingexpodeints-info.tex
\Unterabschnitt{Integrators}{expodeintegrators}

In this section we will describe the steps necessary to create a new \EXPODE integrator.  In distinction to the existing integrators we will use a fictitious integrator called \EXPNEW here.

In the \EXPODEDEV package there are a number of stub functions for the minimally necessary files which need to be created. To create our basic \EXPNEW integrator we type
\begin{Matlabfun}{bash:createNewInt}%
	/path/to/expode \$ cd prototype \\%
	/path/to/expode/prototype \$ ./generateIntegrator.sh expnew "{}new-type"{}
\end{Matlabfun}
Its output will look like this:
\begin{Matlabfun}{bash:createNewIntOutput}%
	Generating integrator expnew, an exponential new-type method. \\%
	[~...~] \\%
	All done. \\%
	Don't forget to edit the ../expnew/expnewOpts.m and ../expnew/expnewSetup.m \\%
	files. Also add the integrator to ../expode/integratorOpts.m.
\end{Matlabfun}
We now have a basic exponential new-type integrator \EXPNEW. To activate the integrator, edit the file \path{/path/to/expode/expode/integratorOpts.m} and add \str{expnew} to the list of available integrators. Let us first look at the files the script has created.

In the \path{/path/to/expode} directory we got the new files \path{expnew.m}, the integrator command, \path{expmsset.m}, the options setter and \path{expnewinfo.m}, the info command. These files are already ready to use, as they only wrap already existing \EXPODE functions. These wrappers will be adapted by \code{generateIntegrator.sh} to give the correct information to their wrapped functions. The script gives us a hint where to start our work. All necessary information for the integration process have to contained in the two functions in \path{/path/to/expode/expnew}, \path{expnewSetup.m} and \path{expnewOpts.m}. The latter one will be discussed in the next section where we handle options.

The \path{expnewSetup.m} has to feed the \EXPODECMD in its startup phase where it will be called the following way:
\begin{Matlabfun}{setupCommandCall}%
	\FUNCTION [ \v{o}, \v{nv}, \v{solverSetup}, \v{solverStep}, \v{order}, \v{errorOrder}, \v{multiStep}, ... \\%
   ~~~~~~~~\v{semilin}, \v{denseOutputGenerator} ] = expnewSetup(\v{o}, \v{nv});
\end{Matlabfun}
The parameters \var{o} and \var{nv} contain the options provided by the user and the numeric values of the list-type options, see sections \refUnterabschnitt{expodeoptions} and \refUnterabschnitt{eD} for details. In case \var{o} or \var{nv} need to be changed by the setup routine, they will be returned again. The arguments \var{solverSetup} and \var{solverStep} will be discussed below in more detail. Next we have \var{order} and \var{errorOrder} which contain the order of the integrator -- needed for the user -- and the order of the error estimator -- needed for the step size selection -- respectively. Futhermore \var{multiStep} counts the number of old time steps to use to approximate the solution at the next step, \var{semilin} indicates, wheter this integrator is built for semilinear problem, \EXPODECMD evaluates the linear part once at the beginning instead of evaluating the Jacobian at each time step. \var{denseOutputGenerator} is a function handle to the dense output generating function. If the integrator has no specific one, this can be empty.

Additionally to these return arguments, the integrator has to set the \code{\var{eD}.\var{int}.\var{jobFunctions}} field. It has to be a cell array of function handles to the scalar, vectorized versions of the matrix functions, that need to be evaluated. The order of the elements has to correspond to the order of the job coefficients that are passed to the matrix function evaluator, see section \refAbschnitt{matFun} for reference.

It can be necessary to access some of the options set by the user and set some data depending on the user's choice. An option's value can be simply accessed via \var{o}.\var{OptName}. Options of type \refOptType{boolean} will be evaluated to \literal{0} for \code{false} and \literal{1} for \code{true}, so you can simply use
\begin{Matlabfun}{expodesetBooleanOption}%
	\IF \var{o}.\var{BooleanOption} \\%
		\c{BooleanOption was set \str{on}} \\%
	\ELSE \ind \\%
		\c{BooleanOption was set \str{off}} \\%
	\END
\end{Matlabfun}
to handle these options. The same works for \refOptType{list} type options. You can use the \var{nv} variable that contains the numeric values corresponding to the list's entries and it is nice to \code{\textcolor{Command}{switch}} the option:
\begin{Matlabfun}{expodesetListOption}%
	\SWITCH \var{o}.\var{List} \\%
		\CASE \var{nv}.\var{List}.\var{value1} \\%
			\c{List was set \str{value1}} \\%
		\CASE \var{nv}.\var{List}.\var{value2} \\%
			\c{List was set \str{value2}} \\%
		\CASE \var{nv}.\var{List}.\var{value3} \\%
			\c{List was set \str{value3}} \\%
	\ENDCASE
\end{Matlabfun}

% Depending on the included options the \EXPODECMD integrator knows, how it should act. It evaluates the Jacobian at each time step or the linear part at the beginning of the integration process respectively and 

We now discuss the missing two return arguments \var{solverSetup} and \var{solverStep} of the setup command, which both are function handles. The first one will be called to initialize the solver step once, after all initialization is complete and the integration process is about to begin, just before the initialization of the matrix function evaluator. This function gets no arguments. For simple, static integrators, where there is no need for such a function, it can be set to \code{@nullfunction}. The \var{solverStep} performs the actual integration step. It has to have the following signature:
\begin{Matlabfun}{solverStepCall}%
	\FUNCTION [ \var{yNew}, \var{normErr}, \var{hOut} ] = solverStep(\var{t}, \var{y}, \var{h}, \var{reuse})
\end{Matlabfun}
Here we obviously have the time \var{t} and phase variable \var{y}, where we start our step with step width \var{h}. The output \var{yNew} is the numerical solution at time $\var{t} + \var{hOut}$. Note that it is possible that the step size had to be reduced from \var{h} to \var{hOut} during the computation of the new step, for instance if a Krylov method was used. The output argument \var{normErr} contains the norm computed by the error estimator, if one is in use, set to \literal{0} otherwise. The \var{reuse} argument is \literal{true} if the preceeding step was rejected due to the error estimator. It should be passed to the matrix function evaluator for all complutations that are only based on the old time step. In case it is \literal{true} some of the previously computed could be reused.

It is generally a good idea to have a look at the existing implementations of solver steps, such as \code{exprb32}, \code{ms} or \code{exp4}, to get a good impression how they have to work. The use of the matrix function evaluators was already explained in section \refAbschnitt{matFun}. Note that it might be possible that the matrix function evaluator requires to reduce the step size. In this case all previously computed data for the next time step that is influenced by the step size have to be recalculated. That is the reason for the \code{\textcolor{Command}{while}}-loops in all the existing solver steps.

You can access old and the current evaluations of the right hand side and possibly of the non-linear part via
\begin{Matlab}{accessF}
  [ \v{F1}, \v{G1} ] = getOldF(\literal{0});
\end{Matlab}
where the argument of \code{getOldF} determines how many steps you need to go to the past, starting with \literal{0} for $F(y_n)$ only.

Now you should have the basic information you need to write \code{expnewStep}. The remainder of this document deals with the definition of user options and the usage of the global \var{eD} variable which are the two crutial parts left to understand to complete the integrator.

%It is best to look at the existing implementation of such solver steps, e.g. \code{exprb32}, \code{ms} or \code{exp4}, to get a good impression how they have to work. The usage of matrix function evaluators was already discussed in section \refAbschnitt{matFun}. Note that if the evaluator decides to reduce the step size all previously computed data that is influenced by the step size have to be recalculated. That is the reason foer the \code{\textcolor{Command}{while}}-loops in all the existing solver steps.

%The script already gives us a hint where to start. The \EXPODE package is built very modular and in such a way, that some of the new files are already useable. The actual integration will be hooked into the \EXPODECMD procedure.

%The latter one will be discussed in the \hrefUnterabschnitt{expodeoptions}{next} section. The first one provides 